
\magnification=1100
\overfullrule0pt

\input prepictex
\input pictex
\input postpictex
\input amssym.def


\def\qed{\hbox{\hskip 1pt\vrule width4pt height 6pt depth1.5pt \hskip 1pt}}

\def\CC{{\Bbb C}}
\def\FF{{\Bbb F}}
\def\RR{{\Bbb R}}
\def\ZZ{{\Bbb Z}}

\def\End{{\rm End}}
\def\id{{\rm id}}
\def\Id{{\rm Id}}
\def\ind{{\rm ind}}
\def\tr{{\rm tr}}
\def\Tr{{\rm Tr}}
\def\Var{{\rm Var}}




\font\smallcaps=cmcsc10
\font\titlefont=cmr10 scaled \magstep1

\font\sectionfont=cmbx10
\font\tinyrm=cmr10 at 8pt


\newcount\sectno
\newcount\subsectno
\newcount\resultno

\def\section #1. #2\par{
\sectno=#1
\resultno=0
\bigskip\noindent{\sectionfont #1.  #2}~\medbreak}

\def\subsection #1\par{\bigskip\noindent{\it  #1} \medbreak}


\def\prop{ \global\advance\resultno by 1
\bigskip\noindent{\bf Proposition \the\sectno.\the\resultno. }\sl}
\def\lemma{ \global\advance\resultno by 1
\bigskip\noindent{\bf Lemma \the\sectno.\the\resultno. }
\sl}
\def\remark{ \global\advance\resultno by 1
\bigskip\noindent{\bf Remark \the\sectno.\the\resultno. }}
\def\example{ \global\advance\resultno by 1
\bigskip\noindent{\bf Example \the\sectno.\the\resultno. }}
\def\cor{ \global\advance\resultno by 1
\bigskip\noindent{\bf Corollary \the\sectno.\the\resultno. }\sl}
\def\thm{ \global\advance\resultno by 1
\bigskip\noindent{\bf Theorem \the\sectno.\the\resultno. }\sl}
\def\defn{ \global\advance\resultno by 1
\bigskip\noindent{\it Definition \the\sectno.\the\resultno. }\slrm}
\def\endthm{\rm\bigskip}

\def\endprop{\rm\bigskip}

\def\pf{\rm\bigskip\noindent{\it Proof. }}
\def\endpf{\qed\hfil\bigskip}


\def\formula{\global\advance\resultno by 1
\eqno{(\the\sectno.\the\resultno)}}
\def\formulano{\global\advance\resultno by 1 (\the\sectno.\the\resultno)}
\def\tableno{\global\advance\resultno by 1
\the\sectno.\the\resultno. }
\def\lformula{\global\advance\resultno by 1
\leqno(\the\sectno.\the\resultno)}

\def\monthname {\ifcase\month\or January\or February\or March\or April\or
May\or June\or
July\or August\or September\or October\or November\or December\fi}

\newcount\mins  \newcount\hours  \hours=\time \mins=\time
\def\now{\divide\hours by60 \multiply\hours by60 \advance\mins by-\hours
     \divide\hours by60         
     \ifnum\hours>12 \advance\hours by-12
       \number\hours:\ifnum\mins<10 0\fi\number\mins\ P.M.\else
       \number\hours:\ifnum\mins<10 0\fi\number\mins\ A.M.\fi}


\nopagenumbers
\def\runningtitle{\smallcaps Metropolis scans and Hecke algebras}
\headline={\ifnum\pageno>1\eoheadline\else\firstheadline\fi}
\def\names{\smallcaps persi diaconis\quad and\quad arun ram}
\def\firstheadline{\noindent Preliminary Draft \hfill  \today}
\def\firstheadline{}
\def\eoheadline{\ifodd\pageno\oddheadline\else\evenheadline\fi}
\def\oddheadline{\tenrm\hfil\runningtitle\hfil\folio}
\def\evenheadline{\tenrm \folio\hfil{\names}\hfil}

\vphantom{$ $}  
\vskip.75truein

\centerline{\titlefont Analysis of systematic scan Metropolis algorithms}
\centerline{\titlefont using Iwahori-Hecke algebra techniques}
\bigskip
\centerline{\rm Persi Diaconis${}^\ast$}
\centerline{Department of Mathematics}
\centerline{Stanford University}
\centerline{Stanford, CA 94305}
\bigskip
\centerline{\rm Arun Ram${}^{\ast\ast}$}
\centerline{Department of Mathematics}
\centerline{University of Wisconsin, Madison}
\centerline{Madison, WI 53706}
\centerline{{\tt ram@math.wisc.edu}}

\footnote{}{\tinyrm ${}^{\ast}$ Research supported in part by National
Science Foundation grant DMS-9504379.}
\footnote{}{\tinyrm ${}^{\ast\ast}$ Research supported in part by National
Science Foundation grant DMS-9622985.}

\bigskip

\noindent{\bf Abstract.}
We give the first analysis of a systematic scan version of the Metropolis
algorithm.  Our examples include generating random elements of a Coxeter group 
with probability determined by the length function.  The analysis is based on
interpreting Metropolis walks in terms of the multiplication in the Iwahori-Hecke algebra.

\section 1. Introduction

When faced with a complex task, is it better to be systematic or proceed by 
making random adjustments?  We study aspects of this problem in the context 
of generating random elements of a finite group.  For example, suppose we want
to fill $n$ empty spaces with zeros and ones such that the probability of 
configuration $x=(x_1,\ldots,x_n)$ is $\theta^{n-|x|}(1-\theta)^{|x|}$
with $|x|$ the number of ones in $x$.  A systematic scan approach works left
to right filling each successive place with a $\theta$-coin toss.  A random 
scan approach picks places at random and a given site may be hit many times
before all sites are hit.  The systematic approach takes order $n$ steps and
the random approach takes order ${1\over 4}n\log n$ steps.

Realistic versions of this toy problem arise in image analysis and Ising
like simulations where one must generate a random array by a Monte Carlo
Markov chain.  Systematic updating and random updating are competing 
algorithms discussed in detail in Section 2.  There are some successful
analyses for random scan algorithms, but the intuitively appealing systematic
scan algorithms have resisted analysis.

Our main results show that the binary problem above is exceptional; for the examples
analyzed in this paper, systematic and random scans converge in about the same 
number of steps.

Let $W$ be a finite Coxeter group generated by simple reflections 
$s_1,s_2,\ldots, s_n$, where $s_i^2=\id$.  For example, $W$ may
be the permutation group $S_{n+1}$ with $s_i=(i,i+1)$.  The length function $\ell(w)$
is the smallest $k$ such that $w=s_{i_1}s_{i_2}\cdots s_{i_k}$.  Fix $0<\theta\le 1$
and define a probability distribution on $W$ by 
$$\pi(w) = {\theta^{-\ell(w)}\over P_{{}_W}(\theta^{-1})},
\qquad\hbox{where}\quad P_{{}_W}(\theta^{-1})=\sum_{w\in W} \theta^{-\ell(w)}
\formula$$
is the normalizing constant.
Thus $\pi(w)$ is smallest when $w=\id$ and, as $\theta\to 1$, $\pi$
tends to the uniform distribution.  These non-uniform distributions
arise in statistical work as Mallows models.  Background and references
are in Section 2e.

A standard Monte Carlo Markov chain algorithm for sampling from $\pi$
is the Metropolis algorithm with a systematic scan.  This algorithm cycles
through the generators in order.  If multiplying by the current generator
increases length this multiplication is made.  If the length decreases,
then the multiplication is made with probability $\theta$ and omitted
with probability $1-\theta$.  One scan uses $s_1, s_2, \ldots, s_{n-1},
s_n, s_n, s_{n-1},\ldots, s_1$, in order.  Define
$$K(w,w') = 
\hbox{the chance that a systematic scan started at $w$ ends in $w'$.}
\formula
$$
Repeated scans of the algorithm are defined by
$$K^\ell(w,w') = \sum_{w''}
K^{\ell-1}(w,w'')K(w'',w'),
\qquad \ell\ge 2.
\formula
$$
In Section 2c and 4a we show that this Markov chain has $\pi$ as unique 
stationary distribution.  

The main results of this paper derive sharp results on rates of convergence
for these walks. 
As an example of what our methods give, we show that order $n$ scans are 
necessary and suffice to reach stationary on the symmetric group 
starting from the identity.  More precisely, we prove

\thm   
Let $S_n$ be the permutation group on $n$ letters.  Fix $\theta$,
$0<\theta\le 1$.  Let $K_1^\ell(w)=K^\ell({\rm id},w)$ be the systematic scan chain 
on $S_n$ defined by (1.2) and (1.3).  
For $\ell=n/2-(\log n)/(\log\theta) + c$ with $c>0$,
$$\Vert K_1^\ell-\pi_\theta\Vert^2_{{}_{TV}} \le
\left(e^{\theta^{2c+1}}-1\right)+n!\theta^{n^2/8-n(\log n)/(\log\theta)+n(c+1/4)}.
\formula
$$
Conversely, if $\ell\le n/4$ then, for fixed $\theta$, $\Vert K_1^\ell-\pi\Vert_{{}_{TV}}$
tends to $1$ as $n\to \infty$.
\endthm
The total variation norm above is defined in Section 2a below.  Note that the 
upper bound in (1.5) tends to zero for $c$ large, so that about $n/2$ 
scans suffice to reach stationarity.  The lower bound shows this is
of the right order for large $n$.

Each scan above uses $2n$ multiplications.  Thus Theorem 1.4 implies that the
systematic scan approach reaches stationarity in $n^2$ operations up to lower 
order terms.  We also conjecture that the random scan approach (see Section 2b)
for this example takes order $n^2$ operations.  Further, in Section 7, we prove 
that the scan based on the sequence 
$$(s_1,s_2,\ldots, s_n, s_n, \ldots, s_1),(s_1,\ldots, s_{n-1}, s_{n-1},\ldots, s_1),
\ldots, (s_1, s_2, s_2, s_1), (s_1, s_1)$$ 
converges in one pass.  Thus, again, up to lower order terms, 
$n^2$ operations suffice to reach stationarity.
These results shows that various different scanning strategies
take the same number of operations to reach stationarity. 

One novel aspect of present arguments is our use of the Iwahori-Hecke
algebra $H$ spanned by the symbols $\{T_w\}_{w\in W}$.  This 
is generated by $T_i=T_{s_i}$, $1\le i\le n$ with the relations
$$T_iT_w = 
\cases{
T_{s_iw}, &if $\ell(s_iw)>\ell(w)$, \cr
q T_{s_iw}+(q-1)T_w, &if $\ell(s_iw)<\ell(w)$. \cr
}$$
We have succeeded in giving an algebraic interpretation of the Markov
chain $K(w,w')$ as multiplication in the Iwahori-Hecke algebra $H$.
From there, knowledge of the center
of $H$  (via a result of Brieskorn-Saito and Deligne)
allows us to explicitly diagonalize $K(w,w')$.
Convergence bounds are given in terms of the eigenvalues and the
generic degrees of representation theory.  Then calculus leads to results 
like Theorem 1.4.

Section 2 collects together probabilistic background and tools.  We explain
Markov chains, the Metropolis algorithm, systematic scans, and relate
the basic Metropolis chain to a natural walk on the 
chambers of a building.  In Section 2e we develop properties
of the measures $\pi$.  Some of these are new even for reflection groups of type A 
(the symmetric group).  These properties will be applied 
to prove lower bounds for walks as in Theorem 1.4.  

Section 3 collects together representation theoretic background and
tools and connects the representation theory to Markov chains.  
Section 4 connects Hecke algebras to the Metropolis algorithm and 
specializes the results from Section 3.  A basic upper bound for
convergence is derived by relating two inner products.  

Sections 5 and 6 derive results for the hypercube and the dihedral groups.
Here we find that both the systematic and random scans converge in about the 
same number of steps -- the differences are only in the lead term constants 
(which are functions of $\theta$).

Section 7 derives results for two different systematic scanning plans 
for the symmetric group.  Though we do not have the space to treat further
examples in this paper, it should be remarked that
the methods of Section 7 should also produce analogous
results for the Weyl groups of type $B_n$ and the imprimitive complex 
reflection groups $G(r,1,n)$. The long and short systematic scans
can be defined in a similar way and the representation theory goes through
without problems (see [Hf], [Ra] and [AK]).  The remaining necessary ingredient 
is an analogue of Lemma 7.2.

\medskip\noindent
{\it Acknowledgement.}  This paper is dedicated to our friend Bill Fulton. 
We are also thankful to Ruth Lawrence for early efforts to help understand
deformed random walks.

\section 2. Probabilistic background

In this section we give background material for Markov chains, the Metropolis 
algorithm, and systematic scans.  In Section 2d we interpret the basic 
walk as a walk on flags and the chambers of a building and in Section 2e
we derive basic properties of the stationary distributions.

\subsection 2a.  Markov chains

Background for Markov chains may be found in any standard probability text
e.g. Feller [Fe, Chapter XV].  For the quantitative theory developed here
see Saloff-Coste [SC] and the references therein.

Let $X$ be a finite set.  A {\it Markov chain} on $X$ is a matrix
$K= \left(K(x,y)\right)_{x,y\in X}$
such that
$$K(x,y)\in [0,1]
\qquad\hbox{and}\qquad
\sum_{y\in X} K(x,y) = 1.$$
The set $X$ is the {\it state space} and $K(x,y)$ gives the probability
of moving from $x$ to $y$ in one step.  Powers of the matrix $K$ give the 
probability of moving from $x$ to $y$ in more steps.  For example,
$$K^2(x,y) = \sum_{z\in X} K(x,z)K(z,y)$$
indicates that to move from $x$ to $y$ in two steps the chain must move to $z$ and
then from $z$ to $y$.  The chain is {\it irreducible and aperiodic} if there is an
$\ell>0$ such that $K^\ell(x,y)>0$ for all $x,y\in X$.
The chain $K$ is {\it reversible} if there is a {\it stationary
distribution} $\pi\colon X\to [0,1]$, 
\enspace $\displaystyle{ \sum_{x\in X} \pi(x)=1}$, \enspace 
such that, for all $x,y\in X$,
$$\pi(x)K(x,y) = \pi(y)K(y,x).$$
For irreducible aperiodic $K$, reversibility implies that, for each $x\in X$,
the real numbers $K^\ell(x,y)$ converge to $\pi(y)$ as $\ell\to \infty$.

The quantitative theory of Markov chains studies the speed of convergence.
The {\it total variation distance} of $K^\ell(x,\cdot)$ to $\pi$ is defined by
$$\Vert K_x^\ell-\pi\Vert_{{}_{TV}}
=\max_{A\subseteq X} \left\vert\sum_{y\in A} K^\ell(x,y)-\pi(y)\right\vert.$$
Using the set $A=\{ y\in X\ |\ K^\ell(x,y)>\pi(y)\}$ it is easily shown that
$$\Vert K_x^\ell-\pi\Vert_{{}_{TV}}
={1\over 2}\sum_{y\in X} |K^\ell(x,y)-\pi(y)|.
\formula
$$
Let $L^2(\pi)$ be the space of functions $f\colon X\to \RR$ with the norm
$$\langle f,g\rangle_{{}_2} = \sum_{x\in X} f(x)g(x)\pi(x).\formula$$
The following lemma provides a relation between the total variation
and the $L^2(\pi)$ norms.  This bound is the primary tool for studying
rates of convergence of Markov chains.

\lemma  Let $f\in L^2(\pi)$.  Then
$
\Vert f\Vert_{{}_{TV}}^2
\le
{1\over 4}\Vert f/\pi \Vert_{{}_2}^2$.
\pf
By the Cauchy-Schwartz inequality,
$$
\Vert f\Vert_{{}_{TV}}^2
= \hbox{$1\over 4$}\left(
\sum_{x\in X} {|f(x)|\over \sqrt{\pi(x)}}
\sqrt{\pi(x)}\right)^2  
\le \hbox{$1\over4$}
\left( \sum_{x\in X} {f(x)^2\over \pi(x)} \right) 
\left( \sum_{x\in X} \pi(x)\right) 
= \hbox{$1\over 4$}
\langle f/\pi,f/\pi\rangle_{{}_2}.
\qquad\hbox{\qed} $$

\subsection 2b.  Systematic scan algorithms

Let $\pi$ be a probability distribution on a finite set $X$ and let
$K_1,K_2,\ldots, K_n$ be Markov chains on $X$ each having stationary
distribution $\pi$.  Then any product 
$K_{i_\ell}K_{i_{\ell-1}}\cdots K_{i_1}$ has stationary distribution
$\pi$ and a choice of an infinite sequence $\{i_\ell\}_{\ell=1}^\infty$
gives a scanning strategy.  A random choice of indices gives a random 
scanning stategy.  If each $K_i$ is reversible
for $\pi$, then $K_1K_2\cdots K_{n-1}K_nK_nK_{n-1}\cdots K_2K_1$ is 
an example of a reversible systematic scanning strategy 
(while $K_1\cdots K_n$ is not necessarily reversible).  

In routine applications of the Metropolis algorithm to image analysis and
Ising-like models the state space has coordinates.  Randomized strategies
choose a coordinate at random and attempt to change it.  Systematic strategies
cycle through the coordinates in various orders.  Fishman [Fi] reviews
the literature on scanning strategies and gives some practical comparison.
The scheme underlying Theorem 1.4 is Fishman's Plan 3.

There has been some rigorous work on rates of convergence for
systematic scans in a related case:  Gaussian distribution of
coordinates with the stochastic updating done by the heat bath algorithm
(also known as Glauber dynamics or the Gibbs sampler).  One fascinating
study by Goodman and Sokal [GS] relates scanning strategies to standard
approaches for solving large linear systems.  They show that the systematic
scan heat bath algorithm is a stochastic analog of the Gauss-Seidel algorithm.
Moreover, they show how previous analyses of Gauss-Seidel give the 
eigenvalues of its stochastic counterpart.  Amit [Am1-2] and Amit and Grenander
[AG] have pushed forward and carried out these ideas to give some
comparison of systematic and randomized sweeps in the Gaussian case.
Their approach uses the fact that the heat bath algorithm is a projection
operator.  In the Gaussian case the problem reduces to the computation
of angles between subspaces of a Hilbert space.  Baronne and Frigessi
[BF] and Roberts and Sahk [RS] are related references.

\subsection 2c. The Metropolis algorithm

The Metropolis algorithm gives a way of changing the stationary distribution
of a given Markov chain into any distribution.  It was invented by Metropolis
et al [MR].  A clear description is in Hammersley and Handscomb [HH] and
a recent survey appears in [DS].

Let $X$ be a finite set.  Let $P(x,y)=P(y,x)$ be a symmetric Markov matrix on $X$
and let $\pi$ be a fixed probability distribution on $X$.  Form a new 
chain by the following recipe:
$$M(x,y) = \cases{
P(x,y), &if $x\ne y$ and $\pi(y)\ge \pi(x)$,\cr
\cr
\displaystyle{P(x,y){\pi(y)\over \pi(x)}}, &if $x\ne y$ and $\pi(y)<\pi(x)$,\cr
\cr
\displaystyle{
P(x,x) +\sum_{\pi(z)<\pi(x)}  P(x,z) \left(1-{\pi(z)\over \pi(x)}\right)
}, &if $x=y$.\cr
}
\formula$$
In words:  
\smallskip

{\narrower{\noindent
Form the Metropolis chain from $x$ by choosing $y$ from $P(x,y)$.
If $\pi(y)\ge \pi(x)$ move to $x$.  If $\pi(y)<\pi(x)$ flip a coin
with chance of heads $\pi(y)/\pi(x)$.  If the coin comes up heads move to $y$.
In all other cases stay at $x$.}

}

\smallskip\noindent
As shown in the references above, the
Metropolis chain is reversible with stationary distribution $\pi$.
It is of practical importance that
the chain $M$ can be run knowing $\pi$ only up to a normalizing constant.
Irreducibility and aperiodicity of $M$ must be checked on a case by case basis.

An example of interest is $X=W$, where $W$ is a finite real reflection group
generated by simple reflections $s_1,s_2,\ldots, s_n$.  Let $P(x,y)$ be 
the Markov chain given by
$$P(x,y) = \cases{
1/n, &if $y=s_ix$ for some $i$, \cr
0, &otherwise. \cr}$$
Here $P(x,y)$ is the usual random walk based on a generating set.  It has
uniform stationary distribution.  
Fix $\theta$, $0<\theta\le 1$ and let $\pi$ be as in (1.1).
The Metropolis construction gives the
Markov chain
$$M(x,y) = \cases{ 
1/n, &if $y=s_ix$ and $\ell(y)>\ell(x)$, \cr
\cr
\theta/n, &if $y=s_ix$ and $\ell(y)<\ell(x)$, \cr
\cr
\displaystyle{{1\over n}\sum_{\ell(s_ix)<\ell(x)} (1-\theta) },
&if $y=x$, \cr
\cr
0, &otherwise, \cr
}
\formula
$$
which has stationary distribution $\pi$.
In Section 4a we demonstrate that this is exactly the chain given by left
multiplication by a uniformly chosen generator $\tilde T_i$ in the
Iwahori-Hecke algebra $H$ with $q=\theta^{-1}$.  Similarly, the systematic scan 
chain of Theorem 1.4 can be interpreted via multiplication in $H$.

Despite its widespread use there has been very limited success in analyzing
the time to stationarity of the Metropolis algorithm. 
In the present paper we carry this out for the random scan Metropolis
algorithm (2.5) on the hypercube (Section 5) and on the dihedral group (Section 6).
Though we have not analyzed the random scan Metropolis algorithm 
on the symmetric group we conjecture that order $n^2$ steps are necessary 
and sufficient to achieve stationarity.  A survey of what is rigorously known 
appears in [DS].

Diaconis and Hanlon [DH] studied the example given by $W=S_n$, the symmetric group,
(so $c(w)=n-[\hbox{\# of cycles in $w$}]$),
with input chain
$$P(x,y) = \cases{
\displaystyle{1\over{n\choose2}}, &if $y=(i,j)x$ for some transposition $(i,j)$, \cr
\cr
0, &otherwise, \cr}
\formula
$$
and stationary distribution $\pi(w) = z\theta^{c(w)}$,
where $z$ is a normalizing constant and
$c(w)$ is the minimum number of transpositions needed to sort $w$.
They showed that all eigenvectors of the resulting Metropolis chain are given by
the coefficients of Jack's symmetric functions (expanded in terms of
the power sum symmetric functions) and they used the corresponding eigenvalues 
to give a complete analysis of the running time.

Similar analyses were carried out in the Ph.D.\ theses of Belsley [Be2] and 
Silver [Si].  They worked in abelian groups with $\pi$ proportional to 
$\theta^{\ell(y)}$ where $\ell$ is the length function with respect to a natural 
set of generators.
In several cases they found that the eigenfunctions were natural deformations of
classical orthogonal polynomials.  Ross and Xu [RX] studied the random scan Metropolis
algorithm on the hypercube using its representation as a random walk on a hypergroup.
It should be emphasized that for other choices of
$\pi$, or in non-group cases, careful analysis of rates of convergence for the Metropolis
algorithm is completely open.  

\subsection 2d.  Some other interpretations of the walks

We have presented Theorem 1.4 in an algorithmic context.  Here we
show how the walk (2.5) arises geometrically on the space of flags and as the
natural nearest neighbor walk on the chambers of a building.
The systematic scan walks have similar interpretations.

Let $\FF_q$ be a finite field. A {\it complete flag}
$F=(0=F_0\subseteq F_1\subseteq
F_2\subseteq \cdots\subseteq F_{n-1}\subseteq F_n=V)$ is a nested
increasing sequence of subspaces of an $n$-dimensional vector space
$V$ over $\FF_q$ with $\dim(F_i)=i$.  A natural random walk on
complete flags may be performed as follows:
\itemitem{} Choose $i$, $1\le i\le n-1$, uniformly.
\itemitem{} Replace $F_i$ by a uniformly chosen subspace $\tilde F_i$
with $F_{i-1}\subseteq \tilde F_i \subseteq F_{i+1}$.
\smallskip\noindent
This walk is symmetric, irreducible and aperiodic.  It thus has the uniform
distribution as its unique stationary measure.  It is instructive to think
of the  ``$q=1$'' case.  Then a flag is a nested increasing
chain $\{i_1\}\subseteq \{i_1,i_2\}\subseteq \cdots\subseteq
\{i_1,i_2,\ldots,i_n\}$ of elements of an $n$ set or, equivalently,
a permutation $(i_1,i_2,\ldots, i_n)$.  In this case the walk is
multiplication by  random pairwise adjacent transpositions.

The space of flags may also be identified as the chambers of a building
of type $A_{n-1}$ and in this formulation the walk is described as follows:
\smallskip
\itemitem{} From a chamber $C$ of the building choose one of the 
adjacent chambers uniformly at random and move there.  
\smallskip\noindent
In the elegant, readable
treatment of buildings by Brown [Bw], he explains that flag space 
may be represented as $G/B$ with $G=GL_n(\FF_q)$ and $B$ the 
subgroup of upper triangular matrices in $GL_n(\FF_q)$.  
Then two flags $g_1B$ and $g_2B$  
differ in the $i$th step as above if and only if $g_1P_i=g_2P_i$ 
where $P_i$ is the parabolic subgroup 
$P_i = B\cup Bs_iB$ (see [Bw, pp. 102-103]).  
Thus, if flags $g_1B$ and $g_2B$ are $i$-adjacent then
$g_2=g_1b$ or $g_2=g_1bs_ib'$ with $b,b'\in B$, and so the 
walk on $G/B$ moves from $gB$ to $gg'B$ with $g'$ uniformly chosen
in $B$ or $Bs_iB$.  In this way, choosing an adjacent chamber of the building
at random produces a $B$-invariant walk on $G/B$.  
Finally, the walk on flags gives rise to a natural walk on the
double coset space $B\backslash G/B$ (described in more detail in
Section 3b).  The double coset space is identifiable with the symmetric
group $S_n$ and the induced Markov chain is given by (2.5) with
$\theta = 1/q$.  A similar story holds for the natural walk on any spherical
building.

\subsection 2e. Properties of the stationary distribution

Suppose that $(X,d)$ is a finite metric space.  A simple way of building
probability models on $X$ is to fix $0< \theta \le 1$ and $x_0\in X$ and define
$$\pi(x) ={ q^{d(x,x_0)}\over P_{{}_X}(q) },
\qquad\hbox{where}\quad
q=\theta^{-1}
\quad\hbox{and}\quad
P_{{}_X}(q)=\sum_{x\in X} q^{d(x,x_0)}
\formula
$$
is a normalizing constant.  When $q=1$ the distribution is uniform.

Models of the form (2.7) were introduced by Mallows [Ma] for the study of
permutations.  He used the length function as a distance, 
$\ell(x^{-1}x_0)=d(x,x_0)$, and estimated $q$ and $x_0$ to match data.
Such Mallows models have had application and development for ranked
and partially ranked data using a variety of metrics [D],
[Cr], [FV], [Mar].  They have
also been used for phylogenetic trees [BHV], classification trees [SB]
and compositions [Det].

One problem in studying Mallows models is that the normalizing
constant $P_{{}_X}(q)$ is uncomputable in general.  In such cases
properties of $\pi$ can be studied by simulation using the Metropolis
algorithm given in Section 2c.

For the examples based on reflection groups the normalizing constants
are known and further there is a simple algorithm available for exact 
generation from $\pi$.  These properties are collected together here.
In each case the properties are illustrated for the permutation group;
some of our results are new for the original Mallows model.  Further,
the properties of $\pi$ (particularly Property 4 below)
are used in proving the lower bounds in Theorem 1.4.

Throughout this section we work with the model 
where the underlying space $X=W$ is a finite Coxeter
group generated by simple reflections, $x_0$ is the
identity element of $W$, and the length function is the distance on $W$.
Thus the model is
$$\pi(w) = {q^{\ell(w)}\over P_{{}_W}(q)},
\qquad\hbox{where}\quad
q=\theta^{-1}
\quad\hbox{and}\quad
P_{{}_W}(q) = \sum_{w\in W} q^{\ell(w)}\formula$$
is the {\it Poincar\'e polynomial} of the group $W$.
It is a classical theorem that the normalizing constant has a simple form:
$$P_{{}_W}(q) = \prod_{i=1}^n {q^{d_i}-1\over q-1}, 
\formula$$
for known integers $d_i$, the {\it degrees} of $W$  (see [Hu, Theorem 3.15]).  For the
symmetric group $S_{n+1}$, $d_i=i+1$ for $1\le i\le n$.  
The Poincar\'e polynomial $P_{{}_W}(q)$ will be used crucially in what follows.

\medskip\noindent
{\bf Property 1.}  {\sl $\pi(w) = \pi(w^{-1})$, since $\ell(w)=\ell(w^{-1})$.}

\smallskip\noindent
This invariance under inversion was first used by Mallows [Ma]
to characterize Mallows models in a larger class of measures as follows:
Suppose $n$ objects are to be ranked by making pairwise comparisons.  Suppose
the true ranking is $1<2<3<\cdots<n$ and a subject ranks objects
$i$ and $j$ correctly with probability $p_{ij}$. Let $Q(w)$ be the chance 
that the comparisons lead to the permutation $w$ given that they are all 
consistent.   
Of course, $Q(w)$ depends on the ${n\choose 2}$ parameters $p_{ij}$.  
Mallows proved that {\it if} $Q(w)=Q(w^{-1})$, then for some real numbers $q$ and $\phi$,
$Q(w) = z q^{\ell(w)}\phi^{r(w)}$ with
$r(w) = \sum iw(i)$ and $z$ a normalizing constant.  He further
showed that the two parameters $q$ and $\phi$
were practically indistinguishable for
large $n$ and suggested setting $\phi=1$, leading to the distribution $\pi(w)$.

\medskip\noindent
{\bf Property 2.}  {\sl Let $J\subseteq \{1,2,\ldots,n\}$ and let $W_J$ be the subgroup
of $W$ generated by the the generators $s_i$ for $i\in J$.  
The group $W_J$ is a {\it parabolic}
subgroup of $W$.  Each coset of $W_J$ in $W$ contains a unique 
coset representative $x_j$ of minimal length [Hu, Prop. 1.10] and
the probability of any such coset is computable via
$$\pi(x_jW_J) = q^{\ell(x_j)}
{P_{{}_{W_J}}(q)\over P_{{}_W}(q)}.\formula$$
}

\smallskip\noindent
As an example suppose that $W$ is the symmetric group $S_n$ generated
by $s_1,s_2,\ldots, s_{n-1}$ where $s_i=(i,i+1)$.  If $J= \{1,2,\ldots, n-2\}$
then $W_J$ is the subgroup of permutations which leave $n$ fixed.  The minimal
length coset representatives $x_j$ for the 
cosets of $W_J$ in $W$ have $j$ in position $n$ and the rest of the entries
in order.  Property 2 says that 
$$\pi\left(\{ w\in S_n\ |\ w(n)=j\}\right) = q^{n-j}{(1-q)\over (1-q^n)}.
\formula$$
Similarly, if $J=\{2,3,\ldots, n-1\}$ Property 2 yields
$$\pi\left(\{w\in S_n\ |\ w(1)=j\}\right) 
= q^{j-1} {(1-q)\over (1-q^n)}.\formula$$
Similar formulas can be derived for the cases where
$J$ consists of the first $j$ or last $j$ elements of $\{1,2,\ldots,n\}$.  

In combination with Property 1, (2.11) also provides a fomula for the probability
of the set of permutations with $j$ in the $n^{\rm th}$ position and
(2.12) gives the probability of the set of permutations with $j$ in the first position.
More generally, one can give formulas for the probability of the set of permutations which
have $1,2,\ldots, j$ in any given relative position.

\medskip\noindent
{\bf Property 3.}  {\sl Let $J_1\supseteq J_2\supseteq \cdots \supseteq J_k=\emptyset$
be a sequence of subsets of $\{1,2,\ldots, n\}$.  Then a sequential algorithm for
generating $w$ in $W$ from $\pi$ is to choose, for each $1\le i\le k-1$,
the minimal length coset representative of a coset of 
$W_{J_{i+1}}$ in $W_{J_i}$, $1\le i\le k-1$ and multiply these together.
If $x_i$ is a minimal length coset representative of a coset 
in $W_{J_i}/W_{J_{i+1}}$ choose $x_i$ with probability
$q^{\ell(x_i)}P_{W_{J_{i+1}}}(q)/P_{W_{J_i}}(q)$.
}

\smallskip
As an example, suppose $W$ is the symmetric group $S_n$ generated by $s_1,s_2,\ldots, s_{n-1}$.
If $J_1\supseteq J_2\supseteq \cdots \supseteq J_{n-1}$ is given by 
$J_i = \{ i,i+1,\ldots, n-1\}$ the algorithm can be realized as the following sequential
procedure: Place symbols down sequentially beginning with $1$.  If symbols
$1,2,\ldots, i-1$ have been placed in some order, place $i$ first with
probability $q^{i-1}(1-q)/(1-q^i)$, second with probability
$q^{i-2}(1-q)/(1-q^i)$, $\ldots$, $i^{\rm th}$ with probability
$(1-q)/(1-q^i)$.  Continuing until all $n$ elements are placed gives
an efficient method of choosing from $\pi$.

An application of this is the following clever algorithm suggested by Pak [Pa]
for generating a uniformly chosen element of $GL_n(\FF_q)$.  Choose
$w\in S_n$ with probability proportional to $q^{\ell(w)}$.  Then form
$b_1wb_2$ with $b_1$ and $b_2$ uniformly chosen in the lower triangular matrices in 
$GL_n(\FF_q)$.  This yields an efficient algorithm
for uniform choice in $GL_n(\FF_q)$.  With obvious modifications this procedure
easily adapts to the other finite groups with a BN pair.

\medskip\noindent
{\bf Property 4.}  {\sl Consider a finite Coxeter group with probability distibution $\pi$
as given in (2.8).  Let $Z$ be the 
random variable given by $Z(w) = \ell(w)$ for $w\in W$.  Then, with $d_i$ as in (2.9),
$$
E_\pi(Z) = {nq\over 1-q}
- \sum_{i=1}^n {d_iq^{d_i}\over 1-q^{d_i}}, 
\qquad\hbox{and}\qquad
\Var_\pi(Z) = {nq\over (1-q)^2}
-\sum_{i=1}^n {d_i^2q^{d_i}\over (1-q^{d_i})^2}. \formula$$
}
\pf
The moment generating function of $Z$ is 
$$
M_Z(t) = E_\pi(e^{tZ}) 
= {1\over P_{{}_W}(q)}\sum_{w\in W} (e^tq)^{\ell(w)} 
= {P_{{}_W}(e^tq)\over P_{{}_W}(q)}
=\prod_{i=1}^n { (1-(e^tq)^{d_i})\over (1-e^tq)}
{(1-q)\over (1-q^{d_i})}. $$
It follows that $Z$ is the sum of independent random variables $Z_1,\ldots,
Z_n$, where
$$M_{Z_i}(t) =  {(1-(e^tq)^{d_i})\over (1-e^tq)}
{(1-q)\over (1-q^{d_i})}.$$
Then
$$E_\pi(Z_i) = {d\over dt}M_{Z_i}(t)\big|_{t=0}
= {q\over 1-q} - {d_iq^{d_i}\over 1-q^{d_i}}$$
and 
$$E_\pi(Z_i^2) = {d^2\over dt^2}M_{Z_i}(t)\big|_{t=0}
={q\over 1-q} 
- {2d_iq^{d_i+1}\over (1-q)(1-q^{d_i})}
+ {2q^2\over (1-q)^2}
- {d_i^2q^{d_i}\over 1-q^{d_i}}.$$
It follows that
$$\Var_\pi(Z_i)={q\over (1-q)^2}
- {d_i^2q^{d_i}\over (1-q^{d_i})^2}.
\qquad\hbox{\qed}$$

We remark that for Coxeter groups of type $A_n$, $B_n$, $D_n$ under 
the probability distribution $\pi$,
$\ell(w)$ has an approximately normal distribution with mean and variance as
in (2.13).  This follows from its representation as a sum of 
independent variables in the proof of Property 4.  
For details, see [D, Ch. 6C, Cor. 1-2].

\section 3. Hecke algebras 

This section introduces Hecke algebras as bi-invariant functions on a group.
We develop the needed Fourier analysis and then specialize to the Iwahori-Hecke
algebras associated to finite Coxeter groups.

\subsection 3a.  Algebras and Fourier analysis

Random walks are traditionally analyzed using Fourier analysis [D].
We find this possible in our examples and here explain the basic tools.

An algebra $H$ over $\CC$ is (split) {\it semisimple} if it is isomorphic
to a direct sum of matrix algebras.  This means that there exists
a finite index set $\hat W$, and positive integers $d_\lambda$, 
$\lambda\in \hat W$, such that 
$$H\cong \bigoplus_{\lambda\in \hat W} M_{d_\lambda}(\CC),$$
where $M_{d_\lambda}(\CC)$ is the algebra of $d_\lambda\times d_\lambda$
matrices with entries in $\CC$.
Fix an isomorphism
$$\phi\colon H\longrightarrow \bigoplus_\lambda M_{d_\lambda}(\CC)
\qquad\hbox{ and define }\qquad
e_{ST}^\lambda = \phi^{-1}(E_{ST}^\lambda),
\qquad \lambda\in \hat W,\enspace 1\le S,T\le d_\lambda,\formula$$
where $E_{ST}^\lambda$ is the matrix in 
$\bigoplus_\lambda M_{d_\lambda}(\CC)$ which has a $1$ in the $(S,T)$
entry of the $\lambda$th block and zeros everywhere else.  The elements
$e_{ST}^\lambda\in H$ are a set of {\it matrix units} for $H$.

The matrix units $\{e_{ST}^\lambda\}$ form a basis of $H$ and we write
$$h=\sum_{\lambda\in \hat W} \sum_{1\le S,T\le d_\lambda}
\rho_{ST}^\lambda(h)e_{ST}^\lambda,
\formula$$
for $h\in H$.  The homomorphisms $\rho^\lambda\colon H\to M_{d_\lambda}(\CC)$
and the linear functionals $\chi_{{}_H}^\lambda\colon H\to \CC$ given by
$$\rho^\lambda(h) = \big( \rho_{ST}^\lambda(h)\big)_{1\le S,T\le d_\lambda}
\qquad\hbox{and}\qquad
\chi_{{}_H}^\lambda(h)= \Tr(\rho^\lambda(h))$$
are the {\it irreducible representations} and the
{\it irreducible characters} of $H$, respectively.  
The homomorphisms $\rho^\lambda$
depend on the choice of $\phi$ but the irreducible characters $\chi_{{}_H}^\lambda$
do not.

A {\it trace} on $H$ is a linear functional $\vec t\colon H\to \CC$ such 
that $\vec t(h_1h_2)=\vec t(h_2h_1)$, for all $h_1,h_2\in H$.
Up to constant multiples there is a unique trace on $M_{d_\lambda}(\CC)$
and this implies that for any trace $\vec t\colon H\to \CC$ on $H$
there are unique $t_\lambda\in \CC$, $\lambda\in \hat W$, such that
$$\vec t = \sum_{\lambda\in \hat W} t_\lambda\chi_{{}_H}^\lambda.
\formula$$
The trace $\vec t$ is {\it nondegenerate} if $t_\lambda\ne 0$ for
all $\lambda\in \hat W$. 
Define a symmetric bilinear form $\langle,\rangle_{{}_H}\colon H\times H\to \CC$
on $H$ by
$$\langle h_1,h_2\rangle_{{}_H} = \vec t(h_1h_2),
\qquad\hbox{for $h_1,h_2\in H$.}$$
The form $\langle,\rangle_{{}_H}$ is nondegenerate if and only if $\vec t$ is
a nondegenerate trace.  

Let $\{T_w\}_{w\in W}$ be a basis of $H$.  The Fourier transform of 
$h=\sum_{w\in W} h_w T_w$ at the representation $\rho$ is 
$$\hat h(\rho) = \sum_{w\in W} h_w\rho(T_w).
\formula$$
The Fourier inversion theorems describe the change of basis matrix
between $\{T_w\}$ and $\{ e_{ST}^\lambda\}$
and recovers $h$ from $\{\hat h(\rho^\lambda)\}_{\lambda\in \hat W}$.

\thm (Fourier inversion and Plancherel)  Let $H$ be a semisimple algebra over
$\CC$ with a nodegenerate trace $\vec t$.  
Let $\{T_w\}_{w\in W}$ be a basis for the algebra $H$.  Let $\{\tilde T_w^*\}_{w\in W}$
be the dual basis with respect to $\langle,\rangle_{{}_H}$, that is
$\langle\tilde T_w^*, T_v\rangle_{{}_H} = \delta_{wv}$.  Then, with notations as
in (3.1-3.4),
$$
e_{ST}^\lambda 
=\sum_{w\in W} t_\lambda\rho_{TS}^\lambda(\tilde T_w^*)T_w \, ,
$$
and, for any $h,h_1,h_2\in H$,
$$h_w = \sum_{\lambda\in \hat W} t_\lambda \Tr(\hat h(\rho^\lambda)\rho^\lambda(\tilde T_w^*)),
\qquad\hbox{for $h\in H$, and}
\formula$$
$$\langle h_1,h_2\rangle_{{}_H} 
= \sum_{\lambda\in \hat W} t_\lambda \Tr\left(\hat h_1(\rho^\lambda)
\hat h_2(\rho^\lambda)\right),
\qquad\hbox{for $h_1,h_2\in H$.}\formula$$
\pf
Since $\vec t$ is nondegenerate, the equation
$\vec t(e_{ST}^\lambda) 
= \sum_{\mu\in \hat W} t_\mu\chi_{{}_H}^\mu( e_{ST}^\lambda ) 
= t_\lambda\delta_{ST}
$
implies that 
$$\left\{ {e_{TS}^\lambda\over t_\lambda} \right\}
\quad\hbox{is the dual basis to}\quad
\left\{ e_{ST}^\lambda\right\}
\qquad\hbox{
with respect to $\langle, \rangle_{{}_H}$.}$$
By (3.2), \enspace
$\displaystyle{\rho_{ST}^\lambda(a) = 
{1\over t_\lambda}\langle a, e_{TS}^\lambda\rangle_{{}_H}. }$
\enspace
Thus
$$
e_{ST}^\lambda 
= \sum_{w\in W} \langle e_{ST}^\lambda,\tilde T_w^*\rangle_{{}_H} T_w
=\sum_{w\in W} t_\lambda\rho_{TS}^\lambda(\tilde T_w^*)T_w.
$$
The equation (3.6) is 
$$h_w = \langle h, \tilde T_w^*\rangle_{{}_H} = \vec t(h\tilde T_w^*)
=\sum_{\lambda\in \hat W} t_\lambda \chi^\lambda_{{}_H}(h\tilde T_w^*)
=\sum_{\lambda\in \hat W} t_\lambda \Tr(\hat h(\rho^\lambda)\rho^\lambda(\tilde T_w^*)),$$
and (3.7) is 
$$\langle h_1,h_2\rangle_{{}_H}
=\vec t(h_1h_2) = \sum_\lambda t_\lambda \chi_{{}_H}^\lambda(h_1h_2)
=\sum_\lambda t_\lambda \Tr(\hat h_1(\rho^\lambda)\hat h_2(\rho^\lambda)).
\qquad\hbox{\qed}$$

\subsection 3b. Coset chains and Hecke algebras

Let $G$ be a finite group, $B$ a subgroup of $G$ and
let $Q$ be a left $B$-invariant probability distribution on $G$.
Right multiplication by random picks from $Q$ induces
a random walk on $G$, 
$$X_0 = x_0, X_1=x_0g_1, X_2=x_0g_1g_2, \ldots,
\formula$$  
which, in turn, induces a process on $B$ cosets, $Y_0, Y_1, Y_2, \ldots$, 
where $Y_i$ is the coset containing $X_i$.
The chain on $G$ produced by right multiplication by random picks from $Q$ is
$\tilde K(x,y)=Q(x^{-1}y)$.  The chance that this chain winds up in an element of $yB$
is $\tilde K(x,yB)=Q(x^{-1}yB)$ and, since $Q$ is left $B$-invariant, 
$\tilde K(x,yB)= \tilde K(xb,yB)$ for any $b\in B$.  This invariance is 
a necessary and sufficient condition for the induced coset process
to be a Markov chain for any starting state $x_0B\in G/B$ (see Theorem 6.32 of [KS]).
If the support of $Q$ is not a coset of a subgroup of $G$ then the chain in (3.8)
is irreducible and aperiodic with uniform stationary distribution.
The resulting coset chain is 
$$K(xB, yB)=Q(x^{-1}yB)
\quad\hbox{with stationary distribution}\quad
\pi(xB) = |B|/|G|.$$

If the probability $Q$ is $B$ bi-invariant then the right process (3.8) on $G$ induces a
process on $B$ double cosets by simply reporting which double coset the element 
$X_i$ is in.  The chance the $G$-chain moves from $x$ to an 
element of $ByB$ in one step is $Q(x^{-1}ByB)$ and since this only depends on the 
double coset of $x$ the induced process is a Markov chain on double cosets
for any starting state $Bx_0B$.  Letting $W$ be a set of coset representatives
for the double cosets of $B$ in $G$, the chain is given by
$$K(w,w')=K(w^{-1}Bw'B),
\qquad\hbox{with stationary distribution}\quad
\pi(w) = |BwB|/|G|,$$
where we view the double coset chain as a Markov chain on the set $W$.

The {\it Hecke algebra} of the pair $(G,B)$ is the subalgebra of 
the group algebra of $G$ consisting of the $B$ bi-invariant functions on $G$,
$$H = \{ f\colon G\to \CC\ |\ f(g)=f(b_1gb_2),
\hbox{\ for $g\in G$, $b_1,b_2\in B$}\}.$$
Background on Hecke algebras may be found in Curtis and Reiner [CR, \S 11D].
If $H$ is commutative then $(G,B)$ is called a Gelfand pair and there
is a well developed probabilistic literature surveyed in [Le1-2]
or [D, Ch. 3F].

Let $W$ be a set of representatives of the double cosets in $B\backslash G/B$.
The functions
$$T_w={1\over |B|}\delta_{BwB},
\qquad w\in W,$$
form a basis of $H$,  
where $\delta_{BwB}$ is the characteristic function of the double coset
$BwB$.  The natural anti-involution on $G$ given by
$g\mapsto g^{-1}$ induces an anti-involution ${}^*\colon H\to H$ on $H$ given by
$T_w\mapsto T_{w^{-1}}$.  The trivial representation of $G$ 
restricts to the {\it index} representation of $H$ given by
$$\rho^{\bf 1}(T_w) = \ind(w),
\qquad\hbox{where}\quad \ind(w) = {|BwB|\over |B|}.
\formula$$
An example to keep in mind is $G= GL_n(\FF_q)$ with $B$ the subgroup of
upper triangular matrices.  Then $W$ is the set of permutation matrices
and $\ind(w)=q^{\ell(w)}$.

Let $L(G/B)=\{f\colon G\to \CC\ |\ f(g)=f(gb) \hbox{\ for $g\in G$, $b\in B$}\}.$
The group $G$ acts on the left of $L(G/B)$ and $H(G/B)$ acts on the
right of $L(G/B)$ by convolution.  The raison d'etre for the 
Hecke algebra is that $H = \End_G(L(G/B))$ and, as $(G,H)$ bimodules,
$$L(G/B) = \sum_{\lambda\in \hat W} G^\lambda\otimes H^\lambda,
\formula$$
where $\lambda$ runs over an index set $\hat W$ of all the irreducible
representations of $H$, $G^\lambda$ is an irreducible $G$-module
and $H^\lambda$ is an irreducible $H$ module, see [CR, (11.25)(ii)].
Centralizers of the action of a finite group, in this case $G$ acting
on $L(G/B)$, are semisimple (our base field is $\CC)$ and therefore 
theory of Section 3a applies to Hecke algebras.
The trace of the action of $H$ on $L(G/B)$ is given by
$$\vec t(T_w) = \cases{
|G|/|B|, &if $w=1$, \cr
0, &otherwise. \cr
}
\formula$$
The decomposition (3.10) yields
$$\vec t = \sum_{\lambda\in \hat W} t_\lambda \chi_{{}_H}^\lambda,
\qquad\hbox{where $t_\lambda = \dim(G^\lambda)$,}\formula$$
and $\chi^\lambda_{{}_H}$ are the irreducible characters of $H$.
Define an inner product on $H$ by
$$\langle h_1,h_2\rangle_{{}_H} = \vec t(h_1 h_2),
\qquad h_1,h_2\in H.$$
The basis 
$$\left\{{T_{w^{-1}}\over\ind(w)}\right\}_{w\in W}
\qquad\hbox{is the dual basis to}\qquad
\left\{ T_w\right\}_{w\in W}\formula$$
with respect to $\langle ,\rangle_{{}_H}$, 
i.e. $\langle (1/\ind(v))T_{v^{-1}},T_w\rangle_{{}_H} = \delta_{vw}$,
for all $v,w\in W$ (see [CR, (11.30)(iii)]).

\subsection 3c.  Iwahori-Hecke algebras

The Hecke algebras associated to finite Chevalley groups $G$ 
and their Borel subgroups $B$
have a remarkable structure theory for their double cosets -- they are indexed
by the elements of a finite Coxeter group $W$.  For example,
in the case of the group $G=GL_n(\FF_q)$ and its Borel subgroup $B$ of 
upper triangular matrices, the group $W$ is the symmetric group.
There are many wonderful references for this material; Brown [Bw], 
Bourbaki [Bou, Ch. IV \S 2 Ex. 22-27], Curtis and Reiner [CR, \S 67-68] 
or Carter [Ca, \S 10.8-10.11].  We develop what we need in this section
and give the relation to probability theory.

Let $W$ be a finite Coxeter group generated by {\it simple
reflections} $s_1,\ldots, s_n$.
These define a length function with $\ell(\id)=0$, $\ell(s_i)=1$, and
$\ell(s_iw)=\ell(w)\pm 1$ for each $w\in W$, $1\le i\le n$.
The {\it Iwahori-Hecke algebra} $H$ corresponding to
$W$ is the vector space with basis $\{T_w\ |\ w\in W\}$ and 
multiplication given by
$$T_iT_w =
\cases{
T_{s_iw}, &if $\ell(s_iw)=\ell(w)+1$,\cr
(q-1)T_w+q T_{s_iw}, &if $\ell(s_iw)=\ell(w)-1$.\cr}\formula$$
where $T_i=T_{s_i}$ for $1\le i\le n$.
When $w=s_i$, 
$T_i^2 = (q-1)T_i+q$
or, equivalently,
$(T_i-q)(T_i+1)=0$.

Let $\hat W$ be an index set for the irreducible representations of $W$.
For each $\lambda\in \hat W$ let 
$\chi^\lambda_{{}_W}$ be the corresponding irreducible character of $W$ and let 
$d_\lambda=\chi^\lambda_{{}_W}(1)$ be the dimension of this representation.
The irreducible representations of the Iwahori-Hecke algebra $H$ 
are in one-to-one correspondence with the
irreducible representations of $W$ in such a way that
if $\chi^\lambda_{{}_H}$ is the character of the irreducible representation
of $H$ indexed by $\lambda\in \hat W$ then
$$\chi^\lambda_{{}_H}(T_w)\big\vert_{q=1} = \chi^\lambda_{{}_W}(w),$$
for all $w\in W$, see [CR, (68.21)].  In particular, the dimension of the irreducible
representation of $H$ indexed by $\lambda$ is $d_\lambda$. 

Define a trace $\vec t\colon H\to \CC$ on $H$ by
$$\vec t(T_w) = \cases{
P_{{}_W}(q), &if $w=1$,\cr
0, &otherwise, \cr}
\qquad\hbox{where}\quad
P_{{}_W}(q) = \sum_{w\in W} q^{\ell(w)}
$$
is the Poincar\'e polynomial of $W$.  Then $\vec t$
is the trace on $H$ given by (3.11) and 
the {\it generic degrees} are the constants $t_\lambda$ defined by 
$$\vec t = \sum_{\lambda\in \hat W} t_\lambda \chi^\lambda_{{}_H}.\formula$$
where $\chi^\lambda_{{}_H}$, $\lambda\in \hat W$, are the irreducible characters
of $H$ (see (3.2)).  
If $\langle , \rangle_{{}_H}\colon H\times H\to \CC$ is the 
inner product on $H$ given by 
$\langle h_1,h_2\rangle_{{}_H}=\vec t(h_1h_2)$, for all $h_1,h_2\in H$, then
$$
\langle T_x, T_{y^{-1}}\rangle_{{}_H} = \delta_{xy}q^{\ell(y)}P_{{}_W}(q),
\qquad\hbox{for all $x,y\in W$,}
\formula$$
see [CR, (68.29)]. 

The ``trivial'' representation $\rho^{\bf 1}$ of the Iwahori-Hecke algebra $H$ is the 
one-dimensional representation corresponding to the trivial representation
of $W$.  For $w\in W$,
$$\rho^{\bf 1}(T_w)=\chi^{\bf 1}_{{}_H}(T_w) = q^{\ell(w)},
\qquad\hbox{and}\quad
\pi = {1\over P_{{}_W}(q)}\sum_{w\in W} T_w$$
is the corresponding minimal central idempotent of $H$ (cf. (3.9)).  Since $t_{\bf 1} = 1$,
$$\vec t(h\pi) = t_{\bf 1}\chi^{\bf 1}(h),
\qquad\hbox{and}\qquad
T_w\pi = q^{\ell(w)}\pi,\formula$$
for all $h\in H$ and $w\in W$, see [CR, (68.23) and (68.28)]. 

Let $\tr$ be the trace of the regular representation of $H$, i.e.
$\tr(h)$ is the trace of the linear transformation obtained from the action
of $h$ on $H$ by left multiplication.  Then
$$\tr = \sum_{\lambda\in \hat W} d_\lambda \chi^\lambda_{{}_H},
\formula$$
where $d_\lambda$ are the dimensions of the irreducible representations of $H$ (see
[CR, (3.37)(iii)]).
\smallskip\noindent
Both traces $\tr$ and $\vec t$ are important in our analysis of Metropolis walks
(see, for example, the proof of Proposition 4.8).

\section 4. Metropolis walks and systematic scans

This section brings together previous results in the form needed to prove
our main theorems.  We show that the various systematic scans are precisely
represented as multiplication in the Iwahori-Hecke algebra.  Then representation
theory yields tractable expressions for the norms involved.

\subsection 4a. Metropolis walks on $W$

Let $W$ be a finite Coxeter group generated by simple reflections $s_1,\ldots, s_n$
and, for each $1\le i\le n$, let $P_i(x,y)$ be the 
Markov chain on $W$ given by
$$P_i(x,y) = \cases{
1, &if $y=s_ix$, \cr
0, &otherwise. \cr}$$ 
Fix $\theta$, $0< \theta\le 1$, and let $\pi$ be as in (1.1)
Then the Metropolis construction produces the Markov chain
$$K_i(x,y) = \cases{
1, &if $y=s_ix$ and $\ell(y)>\ell(x)$, \cr
\theta, &if $y=s_ix$ and $\ell(y)<\ell(x)$, \cr
1-\theta, &if $y=x$. \cr}
\formula$$
The chain $K_i$ can be interpreted as follows:

{\narrower{\noindent
From $w$, try to multiply by $s_i$.  If this increases the length, carry
out the multiplication.  If it decreases the length flip a $\theta$-coin.
If the coin comes out heads carry out the multiplication.  If it comes
up tails the chain stays at $w$.}

}

\smallskip\noindent
Of course, the chain based on a
fixed value of $i$ is not irreducible.  However, any convex linear combination and any
symmetric product of reversible Markov chains with a fixed stationary
distribution is reversible with the same stationary distribution.
If $W$ is the symmetric group then the following chains are reversible for $\pi$:
$$\matrix{
\displaystyle{ {1\over n}\sum_{i=1}^n  K_i }
&\qquad &\hbox{(random scan Metropolis)} \cr
\cr
\displaystyle{ K_1 K_2\cdots K_n K_n\cdots  K_2 K_1 }
&\qquad &\hbox{(short systematic scan)} \cr
\cr
\displaystyle{ (K_1 K_2\cdots  K_n K_n\cdots  K_2 K_1) 
\cdots
(K_1 K_2 K_2 K_1)( K_1 K_1) }
&\qquad &\hbox{(long systematic scan)} \cr
}\formula$$
Note that $K_1 K_2\cdots  K_n$ is an irreducible Markov chain 
with $\pi$ stationary.  However, it is not reversible in general.

The following Theorem (which follows directly from our setup)
shows that many Markov chains on $W$ can be obtained
by left multiplication by elements of $H$ on the basis $\{\tilde T_w\}$.
The remaining results in this subsection provide the necessary tools 
for studying the convergence of these chains by using the representation 
theory of the Iwahori-Hecke algebra $H$.  Though we have chosen to focus 
here on the Iwahori-Hecke algebras related to finite reflection groups $W$,
the results of this section hold in a general Hecke algebra context.

\thm
Let $W$ be a finite Coxeter group and let $H$ be the Iwahori-Hecke algebra
with basis $\{T_w\}_{w\in W}$ as defined in (3.14).
Let 
$$q=\theta^{-1},
\qquad
\tilde T_i = T_i/q,
\qquad\hbox{and}\qquad
\tilde T_w = q^{-\ell(w)}T_w,
\quad\hbox{for $w\in W$.}$$
The Metropolis chain $K_i$ in (4.1) is the same as the matrix of
left multiplication by $\tilde T_i$ with respect
to the basis $\{\tilde T_w\}_{w\in W}$ of $H$:
$$\tilde T_i \tilde T_w = 
\cases{
\tilde T_{s_iw}, &if $\ell(s_iw)>\ell(w)$, \cr
(1-\theta) \tilde T_w + \theta\tilde T_{s_iw}, 
&if $\ell(s_iw)<\ell(w)$. \cr}
\formula$$
\endthm

Identify functions $f\colon X\to \RR$ in $L^2(\pi)$ with
elements of the Iwahori-Hecke algebra $H$ via
$$f = \sum_{x\in W} f(x)\tilde T_x.\formula$$  
The following Proposition shows that we can use the inner product $\langle,\rangle_{{}_H}$
on $H$ (defined in Section 3c) to compute norms in $L^2(\pi)$.  
Coupled with Lemma 2.3 it gives bounds on 
rates of convergence in total variation distance.

\prop  Let $W$ be a finite Coxeter group and let $\pi$ be as in (1.1).
With the identification of $L^2(\pi)$ and the Iwahori-Hecke algebra $H$ given by
(4.5),
$$ \langle f/\pi, g/\pi\rangle_{{}_2}
=\langle f, g^*\rangle_{{}_H},
\qquad\hbox{for all $f,g\in L^2(\pi)$.}$$
where ${}^*\colon H\to H$ is the involutive anti-automorphism of $H$ defined
by $T_w^* = T_{w^{-1}}$.  
\pf
Use the notation \enspace
$\displaystyle{
f=\sum_{x\in W} f(x)\tilde T_x = \sum_{x\in W} f(x) q^{-\ell(x)} T_x
= \sum_{x\in W} f_x T_x. }$  Then, since $\theta=q^{-1}$, (2.2) and (1.1)
give
\enspace
$$
\langle f/\pi, g/\pi\rangle_{{}_2}
=\sum_{x\in W} {f(x)g(x)\over \pi(x)} 
=\sum_{x\in W} {f_x q^{\ell(x)} g_x q^{\ell(x)}\over \theta^{-\ell(x)} }
P_{{}_W}(\theta^{-1}) 
=\sum_{x\in W} f_x g_x q^{\ell(x)} P_{{}_W}(q). 
$$
Thus, by (3.16),
$$\langle f/\pi, g/\pi\rangle_{{}_2}
=\sum_{x\in W} f_x g_x \langle T_x, T_{y^{-1}}\rangle_{{}_H} 
=\sum_{x,y\in W} f_x g_y \langle T_x, T_{y^{-1}}\rangle_{{}_H} 
=\langle f, g^*\rangle_{{}_H}. \qquad\hbox{\qed}
$$

The following lemma shows that the inner product in $L^2(\pi)$, reversibility, 
and the involution $*\colon H\to H$ are simply related.

\lemma  
Let $H$ be the Iwahori-Hecke algebra corresponding to a finite real
reflection group $W$ and let $\pi$ be as in (1.1).  
Let $K$ be a reversible Markov chain on $W$ determined
by left multiplication by an element of $H$ (also denoted by $K$).
The chain $K$ operates on $L^2(\pi)$ by \enspace
$\displaystyle{Kf(x) = \sum_{y\in W} K(x,y)f(y). }$ \enspace
Then the following are equivalent:
\smallskip
\item{(a)}  $K$ is reversible with respect to $\pi$,
\smallskip
\item{(b)}  $K$ is self adjoint with respect to $\langle,\rangle_{{}_2}$,
\smallskip
\item{(c)}  $K=K^*$ in the Iwahori-Hecke algebra $H$,
\smallskip\noindent
where $\langle,\rangle_{{}_2}$ is the norm on $L^2(\pi)$ defined in (2.2)
and ${}^*\colon H\to H$ is the involutive anti-automorphism of $H$ defined
by $T_w^* = T_{w^{-1}}$.  
\pf
If $K$ is reversible then
$$\langle Kf,g\rangle_{{}_2}
=\sum_{x,y\in X} K(x,y)f(y)g(x)\pi(x)
=\sum_{x,y\in X} f(y)K(y,x)g(x)\pi(y)
=\langle f,Kg\rangle_{{}_2},$$
and, conversely, if $K$ is self adjoint then
$$\pi(x)K(x,y)=\langle \delta_x,K\delta_y\rangle
=\langle K\delta_x,\delta_y\rangle = \pi(y)K(y,x),$$
where $\delta_z$ denotes the delta function at $z$, given by
$\delta_z(x)=\delta_{zx}$ (Kronecker delta).
So $K$ is reversible if and only if $K$ is self adjoint.  

If $K$ is self adjoint with respect to $\langle,\rangle_{{}_2}$ then,
by Proposition 4.6,
$$\langle Kf,g^*\rangle_{{}_H}
= \langle Kf/\pi,g/\pi\rangle_{{}_2}
= \langle f/\pi,Kg/\pi\rangle_{{}_2}
= \langle f,(Kg)^*\rangle_{{}_H}
= \langle f,g^*K^*\rangle_{{}_H}.$$
Thus, for all $w\in W$,
$\langle K, T_w\rangle_{{}_H}
=\langle 1, T_w K^*\rangle_{{}_H}
=\langle T_w, K^*\rangle_{{}_H} = \langle K^*, T_w\rangle_{{}_H}.$
So $K=K^*$.
\endpf

The following Proposition is a primary tool for studying rates of convergence of
Markov chains on Iwahori-Hecke algebras.  It bounds the $L^2(\pi)$ norm of 
Lemma 2.3 in terms of characters of the Iwahori-Hecke algebra.
In contrast with the way that random walks are often analyzed (see, for example
[DS2]) the following Proposition also shows that the Markov chain given by $K$ 
can be analyzed without knowing the eigenvalues of $K$ -- it is only necessary 
to compute traces.

\prop  
Let $H$ be the Iwahori-Hecke algebra corresponding to a finite real
reflection group $W$.  Let $K$ be a reversible Markov chain on $W$ with
stationary distribution $\pi$ determined
by left multiplication by an element of $H$ (also denoted by $K$).
Let $K_x^\ell$ denote the Markov chain started at $x$ after $\ell$ steps.
Then
\medskip
\item{(a)} $\displaystyle{
\Vert K_x^\ell/\pi -1\Vert_{{}_2}^2 =
q^{-2\ell(x)}
\sum_{\lambda\ne {\bf 1}} t_\lambda\chi^\lambda_{{}_H}(T_{x^{-1}}K^{2\ell}T_x),
}$
\medskip
\item{(b)} $\displaystyle{
\sum_{x\in W}\pi(x)\Vert K_x^\ell/\pi -1\Vert_{{}_2}^2 =
\sum_{\lambda\ne {\bf 1}} d_\lambda \chi^\lambda_{{}_H}(K^{2\ell}),
}$
\medskip\noindent
where $\chi^\lambda_{{}_H}$ are the irreducible characters,
$t_\lambda$ the generic degrees (3.15), and
$d_\lambda$ the dimensions of the irreducible representations of $H$.
\pf
Equation (3.17) says
$\vec t(h\pi) = t_{\bf 1}\chi_{{}_H}^{\bf 1}(h)$, and 
$\tilde T_w\pi = \pi$, for all $h\in H$ and $w\in W$.
Thus, since $K_x^\ell$ is a probability, Proposition 4.6 gives
$$\eqalign{
1&=\langle K_x^{2\ell}/\pi,1\rangle_{{}_2} 
=\langle K^{2\ell} \tilde T_x, \pi\rangle_{{}_H} 
=\langle K^{2\ell} \tilde T_x, \tilde T_{x^{-1}}\pi\rangle_{{}_H} \cr
&=\vec t(K^{2\ell} \tilde T_x\tilde T_{x^{-1}}\pi) 
=t_{\bf 1}\chi^{\bf 1}_{{}_H}(K^{2\ell} \tilde T_x\tilde T_{x^{-1}}) 
=t_{\bf 1}\chi^{\bf 1}_{{}_H}(\tilde T_{x^{-1}}K^{2\ell} \tilde T_x). \cr
}$$
Then, by Proposition 4.6,
$$
\left\langle K_x^\ell/\pi, K_x^\ell/\pi\right\rangle_{{}_2} 
=\left\langle K^\ell\tilde T_x, (K^\ell\tilde T_x)^*\right\rangle_{{}_H} 
=\left\langle K^\ell\tilde T_x, \tilde T_{x^{-1}}(K^\ell)^*\right\rangle_{{}_H}.
$$
Thus, by (3.12) and Lemma 4.7(c),
$$
\left\langle K_x^\ell/\pi, K_x^\ell/\pi\right\rangle_{{}_2} 
=\left\langle K^\ell\tilde T_x, \tilde T_{x^{-1}}K^\ell\right\rangle_{{}_H} 
=\vec t\left(\tilde T_{x^{-1}}K^{2\ell}\tilde T_x\right) 
=q^{-2\ell(x)}
\sum_{\lambda\in \hat W} t_\lambda\chi^\lambda_{{}_H}(T_{x^{-1}}K^{2\ell}T_x). 
$$
Now (a) follows since
$\left\langle K_x^\ell/\pi-1, K_x^\ell/\pi-1\right\rangle_{{}_2}
=\left\langle K_x^\ell/\pi, K_x^\ell/\pi\right\rangle_{{}_2}-1.
$
Part (b) follows similarly from the following calculation.  Using the 
definition (2.2) of the norm on $L^2(\pi)$,
$$\eqalign{
\sum_{x\in W} \pi(x)\left\langle K_x^\ell/\pi, K_x^\ell/\pi\right\rangle_{{}_2}
&=\sum_{x,y\in W} \pi(x){K^\ell(x,y)K^\ell(x,y)\over \pi(y)} 
=\sum_{x,y\in W} \pi(y){K^\ell(y,x)K^\ell(x,y)\over \pi(y)} \cr
\cr
&=\sum_{y\in W} K^{2\ell}(y,y) 
=\tr(K^{2\ell}) 
= \sum_{\lambda\in \hat W} d_\lambda\chi^\lambda_{{}_H}(K^{2\ell}),  \cr
}$$
where $\tr$ is the trace of the regular representation of $H$ given in (3.18).
\endpf

\subsection 4b.  Systematic scans

One case of Proposition 4.8 which can be analyzed for all finite Coxeter groups $W$
is the case when the Markov chain $K$ is a (generalized) {\it systematic scan}.
This is when $K$ is given by left multiplication by the
element $\tilde T_{w_0}^2$ in the Iwahori-Hecke algebra.  In terms of 
the geometry of the Coxeter group this chain is the Metropolis walk on the 
chambers which tries to move a chamber to its opposite chamber and back again
by successive reflections in the walls of chambers.  Since each step is 
a Metropolis step the chance that the chamber returns to its original
position after one pass is not $1$, but depends on the parameter $\theta$.
In the case when $W$ is the symmetric group this chain is the long systematic scan 
defined in (4.2).

Let $z$ be the sum of all the reflections in $W$.  Then $z$ is a central 
element (since it is a conjugacy class sum) in the group algebra 
of $W$ and thus, by Schur's Lemma, $z$ acts by a constant 
$c_\lambda$ on the irreducible representation of $W$ labeled by 
$\lambda\in \hat W$.  The following well known result shows that 
the element $\tilde T_{w_0}^2$, where $w_0$ is the longest element of $W$, 
is an Iwahori-Hecke algebra analogue of the element $z$.  From the
point of view provided by Theorem 4.3 the following Proposition 
determines the eigenvalues (with their multiplicities) of the systematic
scan Metropolis chain $K$ on $W$.

\prop  Let $z$ be the sum of all the reflections in $W$ and let
$w_0$ be the longest element of $W$.  The element $\tilde T_{w_0}^2$ is in the
center of the Iwahori-Hecke algebra $H$ and
$$\rho^\lambda(\tilde T_{w_0}^2) = q^{c_\lambda-\ell(w_0)}\Id,
\qquad\hbox{where}\quad
c_\lambda = {\chi^\lambda_{{}_W}(z)\over d_\lambda},
$$
$\rho^\lambda$ is the irreducible representation of $H$ indexed by $\lambda$, 
$\chi^\lambda_{{}_W}$ is the irreducible character of $W$ labeled by $\lambda$, and 
$d_\lambda=\chi^\lambda_{{}_W}(1)$ is the dimension of this representation.
\pf
This result is standard, see [Ra, (2.4) and (2.5)], so we only sketch the proof
here.  A result of Brieskorn-Saito [BS]
and Deligne [De] says that $T_{w_0}^2$ is in the center of the
corresponding braid group.  Since the Iwahori-Hecke algebra $H$ is a quotient of
the group algebra of the braid group it follows that 
$T_{w_0}^2$ is in the center of $H$. 
The constant by which it acts on the irreducible
representation labeled by $\lambda$ can be checked as follows.
The element $T_{w_0}^2-q^{\ell(w_0)}$ is divisible by $(q-1)$ and
$$z={T_{w_0}^2-q^{\ell(w_0)}\over q-1}\Big\vert_{q=1}.
\qquad\hbox{Since}\quad
{q^{\ell(w_0)+c_\lambda}-q^{\ell(w_0)}\over q-1}\Big\vert_{q=1} = c_\lambda,$$
and $z$ acts by the constant $c_\lambda$ the element $T_{w_0}^2$ must act by 
the constant 
$q^{\ell(w_0)+c_\lambda}$.  The result of the Proposition
now follows since $\tilde T_{w_0}^2 = q^{-2\ell(w_0)}T_{w_0}^2$. 
An alternative way to obtain the constant
$q^{\ell(w_0)+c_\lambda}$ which $T_{w_0}^2$ acts by is to note that  
$$\det(\rho^\lambda(T_i)) = (-1)^{(d_\lambda-\chi^\lambda_{{}_W}(s_i))/2}
q^{(d_\lambda+\chi^\lambda_{{}_W}(s_i))/2},$$
for all $1\le i\le n$, and this and [Bou, Ch. VI \S 1 Cor. 2] imply that
$$\det(T_{w_0}^2) = q^{2\ell(w_0)d_\lambda+2\chi^\lambda_{{}_W}(z))/2}
=q^{d_\lambda(\ell(w_0)+c_\lambda)}.
\qquad\hbox{\qed}$$ 

Combining Propositions 4.9 and 4.8 gives the following
bounds on the convergence of the systematic scan Metropolis walk on 
a finite Coxeter group $W$.  Explicit analyses of these bounds in examples
are given in Sections 5,6,7.

\thm  Let $H$ be the Iwahori-Hecke algebra corresponding to a finite real
reflection group $W$.  Let $K$ be the systematic scan Metropolis chain on $W$,
i.e. the reversible Markov chain on $W$ with
stationary distribution $\pi$ determined
by left multiplication by the element $\tilde T_{w_0}^2$ of $H$ where $w_0$ is 
the longest element of $W$.
Then
\medskip
\item{(a)} $\displaystyle{
\Vert K_1^\ell/\pi -1\Vert_{{}_2}^2 =
\sum_{\lambda\ne {\bf 1}} t_\lambda d_\lambda 
\theta^{2\ell(\ell(w_0)-c_\lambda)},
}$
\medskip
\item{(b)} $\displaystyle{
\sum_{x\in W}\pi(x)\Vert K_x^\ell/\pi -1\Vert_{{}_2}^2 =
\sum_{\lambda\ne {\bf 1}} d_\lambda^2 \theta^{2\ell(\ell(w_0)-c_\lambda)},
}$
\medskip\noindent
where $\ell(w_0)$ is the length of $w_0$,
$t_\lambda$ are the generic degrees (see (3.15)), 
$d_\lambda$ the dimensions of the irreducible representations of $H$,
and the constants $c_\lambda$ are as given in Proposition 4.9.
\endthm

\section 5. The hypercube

We begin with a simple but instructive example where all details can be carried
out.  We are able to analyze and compare both randomized and systematic scans.
We show that both kinds of scans take order $n\log n$ operations to converge to 
stationarity.  For small values of $\theta$ the systematic scan converges faster
and for $\theta$ close to $1$ the random scan converges faster.

\subsection 5a. Preliminaries

The Coxeter group $W=(\ZZ/2\ZZ)^n$ has generators $s_1,s_2,\ldots, s_n$ and relations
$$s_i^2=1,
\qquad s_is_j=s_js_i, \quad\hbox{for all $1\le i,j\le n$.}$$
The set $X = W = (\ZZ/2\ZZ)^n$ is the space of binary $n$-tuples,
$s_i$ is the vector with $1$ in the $i$th coordinate and
$0$ elsewhere, and
the length function is given by $\ell(x)=|x|=(\hbox{\# of ones in $x$})$.
The longest element of $W$ is $w_0 = s_1s_2\cdots s_n$ and $\ell(w_0)=n$.

The irreducible representations $\rho^\lambda$ of the Iwahori-Hecke algebra of $(\ZZ/2\ZZ)^n$
are all one-dimensional and are indexed by $n$-tuples 
$\lambda=(\lambda_1,\ldots,\lambda_n)$, $\lambda_i\in \{0,1\}$.  Let $|\lambda|
=\lambda_1+\cdots+\lambda_n$.  Then
$$
\rho^\lambda(T_i) = \cases{ q, &if $\lambda_i=0$, \cr
-1, &if $\lambda_i=1$, \cr}
\qquad c_\lambda = n-2|\lambda|,
\qquad\hbox{and}\qquad
t_\lambda = q^{|\lambda|},\formula$$
where $c_\lambda$ and $t_\lambda$ are the constants defined in Proposition 4.9 and
(3.15), respectively.

Fix $0<\theta\le 1$ and let
$$\pi(x) = {q^{\ell(x)}\over P_{{}_W}(q)},
\qquad\hbox{where}\quad
q=\theta^{-1}
\quad\hbox{and}\quad
P_{{}_W}(q)=(1+q)^n\formula$$
is a normalizing constant.
Then $\pi$ is a product measure on $(\ZZ/2\ZZ)^n$ since 
$$\pi(x) 
= \left({q\over1+q}\right)^{\ell(x)}\left({1\over1+q}\right)^{n-\ell(x)}.$$

\subsection 5b. Random scan Metropolis

The random scan Metropolis algorithm proceeds by choosing a coordinate at
random and attempting to change to its opposite mod $2$.  If this results
in a one, the change is made.  If the change in a zero, flip a coin with 
parameter $\theta$.  If the flip comes up heads change the 
chosen coordinate to $0$ and if it comes up tails then the coordinate stays 
at $1$.  The resulting chain is
$$K(x,y) = \cases{
(1/n), &if $\ell(y)=\ell(x)+1$, \cr
\cr
(1/n)\theta, &if $\ell(y)=\ell(x)-1$, \cr
\cr
(\ell(x)/n)(1-\theta), &if $y=x$, \cr
\cr
0, &otherwise. \cr
}
\formula$$
The following theorem shows that order $n\log n$ steps are necessary and
sufficient to reach stationarity.

\thm  Let the random scan Metropolis algorithm on $(\ZZ/2\ZZ)^n$ be
defined by (5.3) with $0<\theta\le 1$.  Then,
for any starting state $x$, and any $\ell$, 
$$\Vert K_x^{\ell}/\pi - 1\Vert_{{}_2}^2
= \sum_{\lambda\ne 0} \theta^{2\lambda\cdot x-|\lambda|}
\left(1-{|\lambda|\over n}(1+\theta)\right)^{2\ell}.
\formula
$$
For $0<\theta<1$ and  $\ell=n(\log n-\log\theta+c)/2(1+\theta)$
with $c>0$,
$$\Vert K_x^\ell -\pi\Vert_{{}_{TV}}^2
\le
\left(e^{e^{-c}}-1\right)+e^{-c/2}.
\formula
$$
The bound in (5.6) is sharp in the sense that if 
$\ell=n(\log n-\log\theta+c)/2(1+\theta)$ then for all $\epsilon>0$
there exists a $c<0$ such that 
$\Vert K_0^\ell - \pi\Vert_{{}_{TV}} > 1-\epsilon$ for
all sufficiently large $n$.
\pf
From the definitions of the irreducible representations of $H$,
$$\eqalign{
\chi_{{}_H}^\lambda\left(\tilde T_{x^{-1}}
((\hbox{$1\over n$})(\tilde T_1+\cdots+\tilde T_n))^{2\ell} 
\tilde T_x\right)
&=n^{-2\ell}q^{-2\ell(x)-2\ell}\chi_{{}_H}^\lambda(T_1+\cdots+T_n)^{2\ell}
\chi_{{}_H}^\lambda(T_x)^2 \cr
&=n^{-2\ell}q^{-2\ell}
\left((n-|\lambda|)q-|\lambda|\right)^{2\ell}q^{2(|x|-\lambda\cdot x)}q^{-2|x|} \cr
&=\left(1-(|\lambda|/n)(1+\theta)\right)^{2\ell}\theta^{2\lambda\cdot x}. \cr
}
$$
The first statement then follows from Theorem 4.10(a) with
the value for $t_\lambda$ given in (5.1).
For the second statement we need to bound the sum on the right hand side 
of (5.5). 
Since $\theta\le 1$, $\theta^{\lambda\cdot x}\le \theta^0$ and 
$$\Vert K_x^\ell/\pi -1\Vert_{{}_2}^2
\le \sum_{j=1}^n {n\choose j}\theta^{-j}\left(1-{j\over n}(1+\theta)\right)^{2\ell}.$$
Break the sum at $n/2$.  For the first half use ${n\choose j}\le n^j/j!$ and
$1-x\le e^{-x}$ to give an upper bound
$$\sum_{j=1}^{n/2} {1\over j!}
\left({n\over\theta}\right)^j e^{-j(1+\theta)2\ell/n}
= \sum_{j=1}^{n/2} {e^{-jc}\over j!} 
\le e^{e^{-c}}-1.$$
For the second half change $j$ to $n-k$ and use the same inequalities to get an 
upper bound
$$\sum_{k=0}^{n/2} {n^k\over k!}
\theta^{k-n}e^{-(n-k)(1+\theta)2\ell/n}
=\sum_{k=0}^{n/2} {1\over k!}
e^{k(\log n+\log\theta + 2\ell(1+\theta)/n)-n\log \theta-2\ell n(1+\theta)/n}
$$
Using $\sum_{k=0}^m {A^k\over k!} \le A^m$ for $A\ge 2$ the bound for the second half
becomes
$$e^{ (n/2)(\log n+\log\theta+2\ell(1+\theta)/n)
-n\log\theta-2\ell n(1+\theta)/n}
=e^{(n/2)(\log n-\log\theta-2\ell(1+\theta)/n)}
=e^{-nc/2}.$$

To show that this upper bound is sharp we use the second moment method.
With respect to the action of $K$ on $L^2(\pi)$ defined in Lemma 4.7
the matrix $K(x,y)$ of (5.3) has an orthonormal basis of eigenfunctions
$$f_\lambda(y) = \theta^{-|\lambda|/2}(-\theta)^{\lambda\cdot y},
\qquad\hbox{with eigenvalues}\quad
1-{|\lambda|\over n}(1+\theta),
\qquad \lambda\in (\ZZ/2\ZZ)^n.\formula$$
Let $e_i\in (\ZZ/2\ZZ)^n$ be the vector with $1$ in the 
$i^{\rm th}$ entry and $0$ elsewhere.
We shall use the test function,
$$T(y) = \sum_i f_{e_i}(y) = \theta^{-1/2}\sum_{i=1}^n (-\theta)^{y_i}
=\theta^{-1/2}(n-|y|(1+\theta))
={n\over \sqrt{\theta}}\left( 1- {|y|(1+\theta)\over n}\right).\formula$$
The expectation and the variance of $T$ with respect to the distribution $\pi$
are 
$$E_\pi(T)=0\qquad\hbox{and}\qquad \Var_\pi(T) = E_\pi(T^2) = n.\formula$$
For $i\ne j$, 
$$f_{e_i}f_{e_j} = f_{e_i+e_j},\qquad\hbox{and}\qquad
f_{e_i}^2
=f_0+{1-\theta\over \sqrt{\theta}}f_{e_i},$$
where the second identity is verified by checking that both sides agree 
when evaluated at each of the two cases: $y$ such that $y_i=1$ and $y$ 
such that $y_i=0$.
We can compute the expectation and variance of $T$ under the 
distribution $K_0^\ell$ as follows:
$$E_{\ell,0}(T) = \sum_y K^\ell(0,y)T(y) = \sum_{i=1}^n (K^\ell f_{e_i})(0)
={n\over \sqrt{\theta}}
\left(1-{1+\theta\over n}\right)^\ell,
\qquad\hbox{and}\quad\qquad \formula$$
$$\eqalign{
\Var_{\ell,0}(T) 
&= E_{\ell,0}(T^2)-E_{\ell,0}(T)^2 \cr
&= E_{\ell,0}\left(\sum_{i=1}^n f_{e_i}^2 +\sum_{i\ne j} f_{e_i}f_{e_j}\right)
-{n^2\over\theta}\left( 1-{1+\theta\over n} \right)^{2\ell} \cr
&=n+{n(1-\theta)\over\theta} \left(1-{1+\theta\over n}\right)^\ell
+{n(n-1)\over\theta} \left(1-{2(1+\theta)\over n}\right)^{\ell}
-{n^2\over\theta} \left( 1-{1+\theta\over n} \right)^{2\ell}. \cr
}$$

We want to use these formulas to show that 
$\ell=(n/2(1+\theta))(\log n-\log \theta + c)$ steps are sharp. 
Fixing $k$ and using $\log(1-x) = -x-x^2/2+O(x^3)$ and 
$e^{-x^2/2} = 1-x^2/2+O(x^4)$,
$$\left(1-{k(1+\theta)\over n}\right)^\ell \sim
e^{\ell(-k(1+\theta)/n-\ell k^2(1+\theta)^2/2n^2)}
\sim \left({\theta e^{-c}\over n}\right)^{k/2}
\left(1-{\ell k^2 (1+\theta)^2\over 2n^2}\right),$$
when $n$ is large. Thus, for 
$\ell=(n/2(1+\theta))(\log n-\log \theta + c)$ and 
$n$ large,
$$\eqalign{
E_{\ell,0}(T) &\sim \sqrt{n} e^{-c/2},
\quad\hbox{and}\quad \cr
\Var_{\ell,0}(T) 
&\sim 
n
+\sqrt{n\over \theta}(1-\theta)e^{-c/2}
+(n-1)e^{-c}\left(1-{\ell 4(1+\theta)^2\over 2n^2}\right)
- {n^2\over \theta}{\theta\over n}e^{-c}\left(1-{2\ell (1+\theta)^2\over 2n^2}\right) \cr
&\sim
n +O_{c,\theta}(\sqrt{n})
+ne^{-c}\left(-{\ell 2(1+\theta)^2\over 2n^2}\right)
-e^{-c} \cr
&= n+O_{c,\theta}(\sqrt{n})+O_{c,\theta}(\log n),\cr
}$$
with the error terms depending on $c$ and $\theta$.  By first choosing
$c$ to be a fixed (large) negative number and then letting $n\to \infty$,
we see that, if $b$ is large, the set $A_b=\{ x\ |\ |T(x)|\le b\sqrt{n}\}$
has probability $1-1/b^2$ under $\pi$ and probability
$O(1/b^2)$ under $K_0^\ell$.  This completes the proof of the last statement.
\endpf

\subsection 5c. Systematic scan Metropolis

\medskip
We turn next to the systematic scan version.  
Order $n\log n$ steps are required here too.
Lest the reader think this contradicts the example which begins this paper, 
we note that the opening example (which
actually corresponds to the heat bath updating setup) replaces each 
coordinate with a freshly chosen pick.  Thus a zero coordinate which is chosen
can remain zero with probability $\theta$.  For the Metropolis version
analyzed here a chosen zero must change to a one.

With notation as in Section 5a, let $N$ be the chain on $\ZZ/2\ZZ$ with
matrix $\pmatrix{0 &1 \cr \theta &1-\theta \cr}$.  On
$(\ZZ/2\ZZ)^n$ define $K_i$ acting as $N$ on the $i$th coordinate.  
Let 
$$K = K_1 K_2 \cdots K_n K_n \cdots K_1.
\formula
$$
This is the systematic scan Metropolis algorithm with stationary distribution
$\pi$.  The following theorem gives bounds on the distance to stationarity.

\thm  Let the systematic scan Metropolis algorithm on $(\ZZ/2\ZZ)^n$ be 
defined by (5.11) with $0<\theta\le 1$.  Then,
for any starting state $x$ and any $\ell$
$$\Vert K_x^\ell/\pi-1\Vert_{{}_2}^2
= \sum_{\lambda\ne 0} \theta^{(4\ell-1)|\lambda|+2\lambda\cdot x}.
\formula
$$
For, $0<\theta<1$ and $\ell = \displaystyle{
{1\over 4}\left({\log n+c\over \log(1/\theta)}+1\right)}$ with $c>0$,
$$\Vert K_x^\ell-\pi\Vert_{{}_{TV}}^2
\le {1\over 4}\left(e^{e^{-c}}-1\right).
\formula
$$
The bound in (5.14) is sharp in the sense that if 
$\ell = \displaystyle{
{1\over 4}\left({\log n+c\over \log(1/\theta)}+1\right)}$
then for all $\epsilon>0$
there exists a $c<0$ such that 
$\Vert K_0^\ell - \pi\Vert_{{}_{TV}} > 1-\epsilon$ for
all sufficiently large $n$.
\pf
From the definitions of the ireducible representations of $H$,
$$\chi_{{}_H}^\lambda\left( \tilde T_{x^{-1}}\tilde T_{w_0}^{2\ell}\tilde T_x\right) 
= \theta^{4\ell |\lambda|+2\lambda\cdot x}.$$
The first statement then follows from Proposition 4.8(a).
Since $\theta \le 1$, $\theta^{2\lambda x}\le \theta^0$, and 
$$\eqalign{
\Vert K_x^\ell/\pi-1\Vert_{{}_2}^2
&=\left(1+\theta^{4\ell-1}\right)^n-1
=\left(1+e^{(4\ell-1)\log\theta}\right)^n-1 \cr
&\le \left(e^{e^{(4\ell-1)\log \theta}}\right)^n - 1
= e^{e^{\log n +(4\ell-1)\log \theta}} - 1,\cr}$$
and the upper bound follows.

The proof of the lower bound is similar to the random scan case and we only
give the variants needed.  As in (5.7) the eigenvectors of the Markov chain 
are 
$$f_\lambda(y) = \theta^{-|\lambda|/2}(-\theta)^{\lambda\cdot y},
\qquad\hbox{but with eigenvalues}\quad
\theta^{2\ell|\lambda|},
\qquad\lambda\in (\ZZ/2\ZZ)^n.$$  
Use the same test statistic $T$ as in (5.8).  The expectation and
the variance of $T$ with respect to the distribution
$\pi$ are the same as in (5.9) and a calculation similar to that in (5.10)
yields that the expectation and variance of 
$T$ under the distribution $M_0^\ell$ are
$$E_{\ell,0}(T) = {n\over \sqrt{\theta}} \theta^{2\ell},
\qquad\hbox{and}\qquad
\Var_{\ell,0}(T) = n\left(1+{1-\theta\over \theta}\theta^{2\ell}
-{1\over \theta}\,\theta^{4\ell}\right).$$
For $\ell = \displaystyle{{1\over 4\log(1/\theta)}(\log n+c)}$,
$$E_{\ell,0}(T) = {\sqrt{n}\over \sqrt{\theta}}e^{-c/2},
\quad \hbox{and}\quad
\Var_{\ell,0}(T) = 
n\left(1+{1-\theta\over\theta\sqrt{n}}e^{-c/2}-{1\over \theta n}e^{-c}\right)
=n(1+O_{c,\theta}(1/\sqrt{n})).$$
This allows the argument to conclude as in the proof of Theorem 5.4.
\endpf

After $\ell$ passes the systematic scan algorithm 
makes $2\ell n$ basic steps.  Thus the results of Theorems 5.4 and 5.12
show that both scanning strategies converge in order $n\log n$ basic steps. 
The following table compares the lead term constants for the two scanning strategies
at various values of $\theta$.
$$\matrix{
\hbox{$\theta$} &\qquad &\hbox{random} &\qquad &\hbox{systematic} \cr
\cr
\hbox{general} &&\displaystyle{n\log(n/\theta)\over 2(1+\theta)} 
&&\displaystyle{n\log n\over 2\log(1/\theta)} \cr
\cr
\displaystyle{1\over 1+\epsilon} &&\displaystyle{n\log n\over 4} 
&&\displaystyle{n\log n\over 2\log(1+\epsilon)} \cr
\cr
{1\over 2} && \displaystyle{n\log 2n\over 3} &&\displaystyle{n\log n\over 2\log 2} \cr
\cr
\epsilon && \displaystyle{n\log(n/\epsilon)\over 2} &&\displaystyle{n\log n\over 2\log(1/\epsilon)} \cr
}$$
We see that the lead term constants make the random scan faster as 
$\theta\to 1$ while the systematic scan is faster as $\theta\to 0$.

\section 6.  The dihedral group.

The hypercube of Section 5 is commutative.  This section treats the simplest
noncommutative example -- the dihedral group $D_{2n}$.  We completely analyze
the convergence of both the randomized and systematic scans.  We find that
both scanning strategies take order $n$ operations to converge to stationarity.

\subsection 6a.  Preliminaries

The dihedral group of order $2n$ is the group $W$ given by generators $s_1, s_2$
and relations
$$ s_1^2=1, \quad s_2^2=1, \quad\hbox{and}\quad
\underbrace{s_1s_2s_1\cdots}_{n\ {\rm factors}} = 
\underbrace{s_2s_1s_2\cdots}_{n\ {\rm factors}}.$$  
This is the group of symmetries of a regular $2n$-gon where
$s_1$ and $s_2$ act by reflection in axes through the center of the $2n$-gon
which form an angle of $2\pi/2n$.  
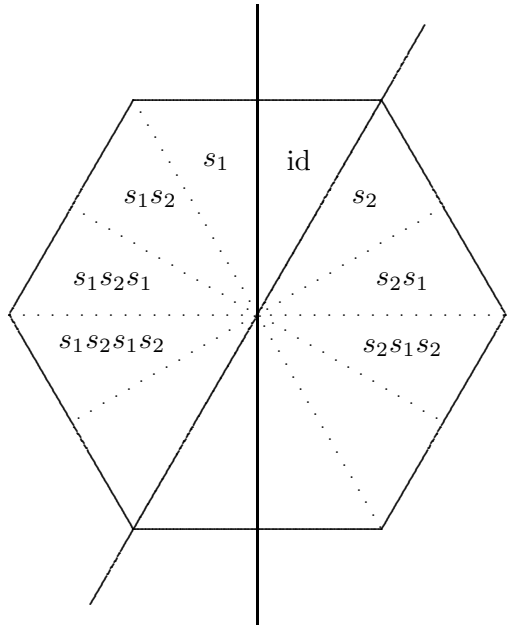
$$
\beginpicture
\setcoordinatesystem units <1cm,1cm>         
\setplotarea x from -5 to 6, y from -4 to 4    
\put{${\rm id}$}[b] at .5 1.75
\put{$s_1$}[b] at -.5 1.75
\put{$s_2$}[b] at 1.3 1.3
\put{$s_2s_1$}[b] at 1.75 .3
\put{$s_1s_2$}[b] at -1.3 1.3
\put{$s_1s_2s_1$}[b] at -1.75 .3
\put{$s_2s_1s_2$}[t] at 1.75 -.3
\put{$s_1s_2s_1s_2$}[t] at -1.76 -.25
\plot -1.5 2.598    1.5 2.598 / 
\plot -1.5 -2.598   1.5 -2.598 /
\plot 1.5 2.598     3 0 /
\plot -1.5 2.598   -3 0 /
\plot -1.5 -2.598  -3 0 /
\plot 1.5 -2.598    3 0 /
\putrule from 0 3.75 to 0 -3.75        
\plot -2.0204 -3.5   2.0204 3.5 /      
\setdots
\plot  -3 0    3 0        /  
\plot -1.5 -2.598  1.5 2.598  /  
\plot -1.5 2.598   1.5 -2.598 /  
\plot -2.25 1.299   2.25 -1.299  /
\plot 2.25 1.299    -2.25 -1.299 /
\endpicture
$$

Fix $0<\theta\le 1$ and consider the distribution on $W$ given by
$$\pi(w) = {q^{\ell(w)}\over P_{{}_W}(q)},
\qquad\hbox{where}\quad
P_{{}_W}(q) = {(q^2-1)\over (q-1)}{(q^n-1)\over (q-1)}
\qquad\hbox{and}\quad q=\theta^{-1}.$$
This measure is largest at the longest element of $W$,
$w_0= s_1s_2s_1\cdots$ ($n$ factors), and $\ell(w_0)=n$.
The walks to be analyzed will all start at the identity.

One may picture the walks described in this section on the $2n$
chambers of an $n$-gon.
Pick one chamber (labeled with identity) and identify
the two internal sides with $s_1,s_2$.  Reflecting the fundamental chamber
around gives each edge and each chamber a label.
The distance $\ell(w)$ is the smallest number of chambers required to walk
from the chamber labeled by $w$ to the identity.  For example, in $D_{12}$, pictured
above $\ell(s_1s_2s_1s_2)=4$.  

The random scan Metropolis walk proceeds from $w$ by choosing one of $s_1, s_2$
with probability $1/2$.  If $\ell(s_iw)>\ell(w)$ the move is accepted.  
If $\ell(s_iw)<\ell(w)$ the move is accepted with probability $\theta$ and rejected
with probability $1-\theta$.

One pass of the systematic scan Metropolis algorithm chooses $n$ generators in the order
$s_1,s_2, s_1, s_2,\ldots$.  Geometrically, starting from the identity, this amounts
to marching around the $n$-gon.  If no rejections are made one complete scan
ends in $w_0$.

Our bounds result in explicit expressions for the convergence of the two walks.
One of these has been carefully analyzed by Belsey [Be, Ch. VI Thm. 2-10].  He showed 

\prop 
For the random scan Metropolis algorithm starting from the identity
$$\Vert K^\ell_1/\pi -1\Vert_{{}_2}^2
\le \theta^{-n}\sqrt{1+\theta\over1-\theta}
\left(1-\hbox{$1\over 2$}(1-\sqrt{\theta})^2\right)^{2\ell}.
\formula$$
For $0<\theta<1$ the right hand side of (6.2) is small for $k$ of order 
$\displaystyle{
{n\log \theta\over 2\log(1-(1/2)(1-\sqrt{\theta})^2)} }$.
\endprop

\noindent
Belsley [Be, Ch. VI Thm. 4-12] further shows that the random scan Metropolis 
algorithm has a total variation cutoff at 
$$\ell = {2n\over 1-\theta} + c\sqrt{n}.$$

For the systematic scan algorithm our results show 
$$\Vert K_1^{\ell}/\pi - 1\Vert_{{}_2}^2
= \theta^{(4\ell-1)n}+\theta^{(2\ell-1)n}\left(
{\theta^2-1\over \theta-1}\cdot{\theta^n-1\over \theta-1} - 1
\right)
-\theta^{2\ell n},$$
and so, when $\ell=1$,
$$\Vert K_1^{\ell}/\pi - 1\Vert_{{}_2}^2
\le \theta^n\left({\theta^2-1\over 1-\theta}\cdot{-1\over \theta-1} - 1\right)
= {2\theta^{n+1}\over 1-\theta}.
\formula$$
Thus, for large $n$ and fixed $0<\theta\le 1$, a single scan suffices to achieve
randomness.

A comparison of the results in (6.2) and (6.3) shows a mild advantage for 
systematic scans.  The effect is most
pronounced as $\theta$ approaches $1$.

\subsection 6b. Representation theory

References for the statements in this paragraph are [KSo] and [CR, \S 67C].
The two dimensional irreducible representations of the Iwahori-Hecke
algebra $H$ are indexed by $0<\lambda<n/2$ and are given explicitly by
$$\rho^\lambda(T_1) = 
q^{1/2}\pmatrix{-d\xi &(1+ad)\xi \cr
\xi^{-1} &-a\xi^{-1} \cr},
\qquad\hbox{and}\qquad
\rho^\lambda(T_2) = 
q^{1/2}\pmatrix{a &1+ad \cr
1 &d \cr},\formula$$
where 
$$\xi=e^{2\pi i\lambda/n}, 
\qquad\hbox{and the equations}\quad
a+d=q^{1/2}-q^{-1/2}
\quad\hbox{and}\quad
-a\xi^{-1}+d\xi=q^{1/2}-q^{-1/2}$$
determine $a$ and $d$. 
For these representations
$$c_\lambda=0
\qquad\hbox{and}\qquad
t_\lambda = {1\over n}{(q^n-1)\over (q-1)}
{(1-\xi)(1-\xi^{-1})\over (q-\xi)(q-\xi^{-1})}q(q+1).
\formula
$$
If $n$ is odd there are two one-dimensional representations of $H$,
the ``trivial'' representation and the ``sign'' representation,
$$\eqalign{
\rho^{\bf 1}(T_1)&=\rho^{\bf 1}(T_2) = q,\cr
\cr
\rho^{\rm sgn}(T_1)&=\rho^{\rm sgn}(T_2) = -1,\cr}
\qquad\eqalign{
\hbox{with}\cr
\cr
\hbox{with}\cr} \qquad
\eqalign{
c_{\bf 1} &= n, \cr
\cr
c_{\rm sgn} &= -n, \cr}
\quad\eqalign{
\hbox{and}\cr
\cr
\hbox{and}\cr} \quad
\eqalign{
t_{\bf 1} &= 1, \cr
\cr
t_{\rm sgn} &= q^n. \cr}
$$
If $n$ is even there are two {\it additional} one dimensional representations
of $H$,
$$\eqalign{
\rho^+(T_1) &= q, \cr
\rho^+(T_2) &= -1, \cr}
\qquad\hbox{and}\qquad
\eqalign{
\rho^-(T_1) &= -1, \cr
\rho^-(T_2) &= q, \cr}
$$
with
$$c_+=c_-=0,
\qquad\hbox{and}\qquad
t_+=t_-={2q(q^n-1)\over n(q^2-1)}.$$

The next Proposition computes the traces on the Iwahori-Hecke algebra which are
needed in order to use Proposition 4.8 to give upper bounds for the convergence
of the random and systematic scan Metropolis walks on the dihedral group.

\prop Let $W$ be the dihedral group of order $2n$ and let $w_0$ be the
longest element of $W$.  Let $\theta = q^{-1}$ and $\tilde T_w = q^{-\ell(w)}T_w$
as in Theorem 4.3.  If $\chi^\lambda_{{}_H}$ is the character of a 
two dimensional representation of the Iwahori-Hecke algebra $H$ then
$$\eqalign{
\chi_{{}_H}^\lambda\left((\tilde T_{w_0}^2)^{2\ell}\right) 
&= 2\theta^{2\ell n}, \qquad\hbox{and} \cr
\chi_{{}_H}^\lambda\left(((1/2)(\tilde T_1+\tilde T_2))^{2\ell}\right) 
&= 
\left( {\theta-2\cos(\pi\lambda/n)\theta^{1/2}-1\over 2}\right)^{2\ell}
+
\left( {\theta+2\cos(\pi\lambda/n)\theta^{1/2}-1\over 2}\right)^{2\ell},
\cr}
$$
For the one dimensional irreducible characters of $H$,
$$\eqalign{
&\chi_{{}_H}^{\bf 1}\left((\tilde T_{w_0}^2)^{2\ell}\right) = 1, \cr
&\chi_{{}_H}^{\bf 1}\left(((1/2)(\tilde T_1+\tilde T_2)^{2\ell}\right) = 1, \cr
}
\qquad\qquad
\eqalign{
&\chi_{{}_H}^{\rm sgn}\left((\tilde T_{w_0}^2)^{2\ell}\right) = \theta^{4\ell n}, \cr
&\chi_{{}_H}^{\rm sgn}\left(((1/2)(\tilde T_1+\tilde T_2))^{2\ell}\right) = \theta^{2\ell}, \cr
}
$$
and, if $n$ is even,
$$\eqalign{
\chi_{{}_H}^+\left((\tilde T_{w_0}^2)^{2\ell}\right) &= 
\chi_{{}_H}^-\left((\tilde T_{w_0}^2)^{2\ell}\right) = \theta^{2\ell n},
\qquad\hbox{and} \cr
\chi_{{}_H}^+\left(((1/2)(\tilde T_1+\tilde T_2))^{2\ell}\right) &= 
\chi_{{}_H}^-\left(((1/2)(\tilde T_1+\tilde T_2))^{2\ell}\right) 
= \left({\theta-1\over 2}\right)^{2\ell}. \cr}
$$
\pf
The formulas for $\chi_{{}_H}^\lambda((\tilde T_{w_0}^2)^{2\ell})$ follow from Proposition 4.9.
The results for the one dimensional representations are easy consequences
of the definitions of the representations $\rho^\lambda$.  Let
us compute $\chi_{{}_H}^\lambda((T_1+T_2)^m)$ for the
two dimensional representations.
From the definition of the representations $\rho^\lambda$ in 6.4
$$\Tr\left(\rho^\lambda(T_1+T_2)\right)
=2q^{1/2}(q^{1/2}-q^{-1/2})$$
and
$$\eqalign{
\det(\rho^\lambda(T_1+T_2))
&=q\left((a-d\xi)(d-a\xi^{-1})+(1+\xi^{-1})(1+ad)(1+\xi)\right) \cr
&=q\left((a+d)(-a\xi^{-1}-d\xi)-(1+\xi^{-1})(1+\xi)\right) \cr
&=q\left((q^{1/2}-q^{-1/2})^2-(\xi^{1/2}+\xi^{-1/2})^2\right). \cr
}$$
where $\xi=e^{2\pi i\lambda/n}$.
Thus the characteristic polynomial of $\rho^\lambda(T_1+T_2)$ is
$$\eqalign{
t^2-&2q^{1/2}(q^{1/2}-q^{-1/2})t
+q\left( (q^{1/2}-q^{-1/2})^2-(\xi^{1/2}+\xi^{-1/2})^2 \right) \cr
&= \left( t-q^{1/2}((q^{1/2}-q^{-1/2})+(\xi^{1/2}+\xi^{-1/2}))\right)
\left( t-q^{1/2}((q^{1/2}-q^{-1/2})-(\xi^{1/2}+\xi^{-1/2}))\right). \cr
}$$
This determines the eigenvalues of $\rho^\lambda(T_1+T_2)$ and
the results follow after substitution of $q=\theta^{-1}$ and $\tilde T_i = q^{-1}T_i$.
\endpf

The next two theorems give explicit expressions for the norms 
$\Vert K_1^\ell/\pi-1\Vert_{{}_2}^2$ for the random and systematic scan Metropolis walks
on the dihedral group.   Bounds on rates of convergence based on these
expressions appear in the summary of results in Section 6a above.

\thm  Let $K$ be the systematic scan Metropolis algorithm on the dihedral group
defined by multiplication by the element 
$\tilde T_{w_0}^2 = \tilde T_1\tilde T_2\tilde T_1\cdots$ ($2n$ factors)
in the Iwahori-Hecke algebra.  Then
\smallskip
\item{(a)} $\displaystyle{
\Vert K_1^\ell/\pi-1\Vert_{{}_2}^2 
= \theta^{(4\ell-1)n}+\theta^{(2\ell-1)n}\left(
{\theta^2-1\over \theta-1}\cdot{\theta^n-1\over \theta-1} - 1
\right)
-\theta^{2\ell n}, }$ \enspace and
\smallskip
\item{(b)} $\displaystyle{
\sum_{x\in W} \pi(x)\Vert K_x^\ell/\pi-1\Vert_{{}_2}^2
= \theta^{4\ell n}+(2n-2)\theta^{2\ell n}. }$
\pf
From (3.15) we have that
$$\sum_{\lambda} t_\lambda d_\lambda = P_{{}_W}(q) = (1+q)\left({q^n-1\over q-1}\right)
=\theta^{-n}(\theta+1)\left({\theta^n-1\over \theta-1}\right).$$
From Theorem 4.10(a)
$$\eqalign{
\Vert K_1^\ell/\pi-1\Vert_{{}_2}^2
&=
\sum_{\lambda\ne {\bf 1}} t_\lambda d_\lambda \theta^{2\ell(\ell(w_0)-c_\lambda)} 
=\left(\sum_{\lambda} t_\lambda d_\lambda \theta^{2\ell n}\right)
-t_{\bf 1}d_{\bf 1}\theta^{2\ell n}
+t_{\rm sgn} d_{\rm sgn}(\theta^{4\ell n}-\theta^{2\ell n}) \cr
&=\theta^{2\ell n}P_{{}_W}(q)-\theta^{2\ell n}+\theta^{-n}(\theta^{4\ell n}-\theta^{2\ell n}) \cr
&= \theta^{(4\ell-1)n}+\theta^{(2\ell-1)}\left(
{\theta^2-1\over \theta-1}\cdot{\theta^n-1\over \theta-1} - 1
\right)
-\theta^{2\ell n}. \cr
}$$
Similarly, use the fact that $\sum_{\lambda\in \hat W} d_\lambda^2 = 2n$ to obtain
$$\eqalign{
\sum_{x\in W} \pi(x)\Vert K_x^\ell/\pi-1\Vert_{{}_2}^2
&=
\sum_{\lambda\ne {\bf 1}} d_\lambda^2 \theta^{2\ell(\ell(w_0)-c_\lambda)} 
=\left(\sum_{\lambda} d_\lambda^2 \theta^{2\ell n}\right)
-d_{\bf 1}^2\theta^{2\ell n}
+d_{\rm sgn}^2(\theta^{4\ell n}-\theta^{2\ell n}) \cr
&=\theta^{2\ell n}2n-\theta^{2\ell n}+\theta^{4\ell n}-\theta^{2\ell n} 
= \theta^{4\ell n}+(2n-2)\theta^{2\ell n}.
\qquad\hbox{\qed} \cr
}$$

\thm  Let $K$ be the random scan Metropolis algorithm on the dihedral group
defined by multiplication by the element $(1/2)(\tilde T_1+\tilde T_2)$
in the Iwahori-Hecke algebra.  Then
$$
\eqalign{
\Vert K_1^\ell/\pi-1\Vert_{{}_2}^2 
=\theta^{2\ell-n}
+{\theta^{1-n}\over n} &\left({\theta^2-1\over \theta-1}\right)
\left({\theta^n-1\over \theta-1}\right) \cr
&\times
\sum_{0<\lambda<n} 
{2-2\cos(2\pi\lambda/n)\over 
(\theta^2-2\cos(2\pi\lambda/n)\theta+1)}
\left({\theta+2\cos(\pi\lambda/n)\theta^{1/2}-1\over 2}\right)^{2\ell},
\cr}
$$
and
$$\sum_{x\in W} \pi(x)\Vert K_x^\ell/\pi-1\Vert_{{}_2}^2
=\theta^{2\ell}
+\sum_{0<\lambda<n} 
2\left({\theta+2\cos(\pi\lambda/n)\theta^{1/2}-1\over 2}\right)^{2\ell}. 
$$
\pf
Substituting $q=\theta^{-1}$ in the formulas for $t_\lambda$ gives
$t_{\rm sgn}=\theta^{-n}$,
$$t_+=t_- = {2\over n}\theta^{1-n}\left({\theta^n-1\over \theta^2-1}\right) 
\qquad\hbox{and}\qquad
t_\lambda=
{\theta^{1-n}(\theta+1)\over n}
{2-2\cos(2\pi\lambda/n)\over (\theta^2-2\cos(2\pi\lambda/n)\theta+1)}
\left({\theta^n-1\over \theta-1}\right),
$$
for $0<\lambda<n/2$.  Thus, if $K=(1/2)(\tilde T_1+\tilde T_2)$,
Theorem 4.10a gives
$$\eqalign{
\Vert K_1^\ell/\pi-1\Vert_{{}_2}^2
&=
\sum_{\lambda\ne {\bf 1}} t_\lambda \chi^\lambda(K^{2\ell})
=
t_+\left({1-\theta\over 2}\right)^{2\ell}+t_-\left({\theta-1\over 2}\right)^{2\ell}
+t_{\rm sgn}\theta^{2\ell} 
+\sum_{0<\lambda<n/2} t_\lambda\chi^\lambda(K^{2\ell}) \cr
&=
\theta^{1-n}{4\over n}\left({\theta^n-1\over \theta^2-1}\right)
\left({\theta-1\over 2}\right)^{2\ell}+\theta^{2\ell-n} 
+\sum_{0<\lambda<n/2} t_\lambda\chi^\lambda(K^{2\ell}), \cr
}$$
where the first term appears only if $n$ is even.  The result now follows
from Proposition 6.6.
Similarly,
$$\eqalign{
\sum_{x\in W} \pi(x)\Vert K_x^\ell/\pi-1\Vert_{{}_2}^2
&=
\sum_{\lambda\ne {\bf 1}} d_\lambda \chi^\lambda(K^{2\ell}) 
=
2\left({1-\theta\over 2}\right)^{2\ell}
+\theta^{2\ell} +\sum_{0<\lambda<n/2} 2\chi^\lambda(K^{2\ell}) \cr
&=\theta^{2\ell}
+\sum_{0<\lambda<n} 
2\left({\theta+2\cos(\pi\lambda/n)\theta^{1/2}-1\over 2}\right)^{2\ell}. 
\quad\hbox{\qed}\cr}$$

\section 7.  The symmetric group

This section proves Theorem 1.4 and a similar result for a different
scanning strategy.  The results show that both scanning strategies
require $n^2$ operations up to lead term constants.

\subsection Preliminaries

The symmetric group $S_n$ is generated by the simple transpositions
$s_i=(i,i+1)$, $1\le i\le n-1$, and the longest element of $S_n$ is the
reversal permutation
$$w_0=\pmatrix{1 &2 &\cdots &n-1 &n \cr n &n-1 &\cdots &2 &1\cr}
\qquad\hbox{with}\qquad \ell(w_0)={n\choose 2}.$$
The book of Fulton [Fu] provides a review of the representation
theory of $S_n$ and we will adopt the conventions for tableaux used there.
The irreducible representations of the Iwahori-Hecke algebra $H$ are indexed by
partitions $\lambda$ with $n$ boxes.
$$
\beginpicture
\setcoordinatesystem units <0.5cm,0.5cm>         
\setplotarea x from 0 to 4, y from 0 to 3    
\linethickness=0.5pt                          
\put{$\lambda=$} at -1 3.5
\putrule from 0 6 to 5 6          %
\putrule from 0 5 to 5 5          
\putrule from 0 4 to 5 4          %
\putrule from 0 3 to 3 3          %
\putrule from 0 2 to 3 2          %
\putrule from 0 1 to 1 1          %
\putrule from 0 0 to 1 0          %

\putrule from 0 0 to 0 6        %
\putrule from 1 0 to 1 6        %
\putrule from 2 2 to 2 6        %
\putrule from 3 2 to 3 6        
\putrule from 4 4 to 4 6        %
\putrule from 5 4 to 5 6        %
\endpicture
$$
Number the rows and columns of $\lambda$ as for matrices.
If $\lambda_i$ and $\lambda_j'$ denote the length of the $i^{\rm th}$ row and
$j^{\rm th}$ column of $\lambda$ respectively, the {\it content} and the
{\it hook length} of a box $b$ in position $(i,j)$ of $\lambda$ are 
$$c(b)=j-i
\qquad\hbox{and}\qquad h_b=\lambda_i-i+\lambda_j'-j+1,$$
respectively.
Set 
$$\displaystyle{n(\lambda)=\sum_{i=1}^{\ell(\lambda)} (i-1)\lambda_i},
\qquad\hbox{and let}\quad
[k]_q = {q^k-1\over q-1}
\quad\hbox{and}\quad
[k]_q!=[k]_q[k-1]_q\cdots [2]_q[1]_q,$$
for each positive integer $k$.
With these notations the dimensions $d_\lambda$ of the irreducible representations
of $H$, the generic degrees $t_\lambda$ defined in (3.15), and the constants
$c_\lambda$ defined in Proposition 4.9, are given by
$$d_\lambda = {n!\over \displaystyle{\prod_{b\in \lambda} h_b} },
\qquad
c_\lambda = \sum_{b\in \lambda} c(b)
\qquad\hbox{and}\qquad
t_\lambda=q^{n(\lambda)}{[n]_q!\over \displaystyle{\prod_{b\in\lambda} [h_b]_q} }\ ,
\formula$$
respectively (see [Fu, \S 7.2 Prop. 2], [Mac, I \S7 Ex. 7 and \S1 Ex. 3], [Hf, 3.4.14]
and [Mac, IV (6.7)]).
The dimension $d_\lambda$ is also equal to the number of standard tableaux of
shape $\lambda$, i.e. fillings of the boxes of $\lambda$ with $1,2,\ldots, n$ such 
that the rows are increasing left to right and the columns are increasing top to bottom
(see [Fu, p.53]).

The next lemma provides bounds on the constants in (7.1) which will be
useful for proving bounds on the convergence times of the systematic
scan Metropolis walks that we analyze here.
The bounds on $c_\lambda$ given in part (c) are essentially those given by
Diaconis and Shahshahani, see [D, 3D Lemma 2].

\lemma
For each partition $\lambda$ let $t_\lambda$, $c_\lambda$ and $d_\lambda$ 
be as defined in (7.1).  
\medskip
\itemitem{(a)}\ \ When $\theta=1/q$ and $0<\theta\le 1$, \enspace
$\displaystyle{
t_\lambda \le \theta^{{\lambda_1\choose 2}-{n\choose 2}}d_\lambda,
}$
\medskip
\itemitem{(b)}\ \ 
$\displaystyle{\sum_{\lambda\vdash n} d_\lambda^2 = n!},$
\qquad and\qquad
$\displaystyle{
\sum_{\lambda\atop \lambda_1=n-j} d_\lambda^2 \le {n^{2j}\over j!},
}$
\qquad
for each $1\le j\le n$.
\medskip
\itemitem{(c)}\ \ 
$\displaystyle{
c_\lambda \le 
\cases{
\displaystyle{ {\lambda_1\choose 2}+\hbox{$1\over2$}(n-\lambda_1)(n-\lambda_1-3), }
&if $\lambda_1\ge n/2$, \cr
\cr
n^2/4-n, &if $\lambda_1\le n/2$. \cr
} }$
\pf 
Set $\theta = 1/q$ and use [Mac, I \S1 Ex. 2] and [Mac, III \S6 Ex. 2], to get
$$
t_\lambda = \theta^{-n(\lambda)-{n\choose 2}-n+\left(\sum h_b\right)}
{[n]_\theta!\over \displaystyle{\prod_{b\in\lambda} [h_b]_\theta} } 
= \theta^{n(\lambda')-{n\choose 2}}
{[n]_\theta!\over \displaystyle{\prod_{b\in\lambda} [h_b]_\theta} } 
= \theta^{-{n\choose 2}}
\sum_Q \theta^{r(Q)},$$
where the sum is over all standard tableaux $Q$ of shape
$\lambda$ and $r(Q)$ is the sum of $i$ such that $i+1$ is to the right of
$i$ in $Q$.  Thus $t_\lambda$ is a sum of 
$d_\lambda$ monomials where the lowest degree term has degree $n(\lambda')-{n\choose2}
\ge {\lambda_1\choose2}-{n\choose2}$.  Part (a) follows.  

We can bound the number of standard tableaux $Q$ of shape $\lambda$
with $\lambda_1=n-j$ by noting that
there are ${n\choose j}$ ways of picking the elements not in the
first row of $Q$ and at most $\sqrt{j!}$ ways of arranging
these to complete a standard tableau.  Thus
$$\sum_{\lambda\atop\lambda_1 = n-j} d_\lambda^2
\le \left(\sum_{\lambda_1=n-j} d_\lambda\right)^2
\le \left( {n\choose j}\sqrt{j!}\right)^2
\le {n^{2j}\over (j!)^2}j! = {n^{2j}\over j!}.
$$

The inequalities in (c) are direct consequences of 
$$c_\lambda \le 
c_{(\lambda_1,n-\lambda_1)},\enspace\hbox{if $\lambda\ge n/2$,\qquad and}\qquad
c_\lambda\le c_{(n/2,n/2)},\enspace\hbox{if $\lambda_1\le n/2$.}
\qquad\hbox{\qed}$$

\subsection Long systematic scan

As in Section 4a we fix $0<\theta\le 1$ and consider the Markov chain
$$K_i(x,y) = \cases{
1, &if $y=s_ix$ and $\ell(y)>\ell(x)$, \cr
\theta, &if $y=s_ix$ and $\ell(y)<\ell(x)$, \cr
1-\theta, &if $y=x$. \cr}
\formula
$$
which is produced by applying the Metropolis construction to the base chain
$$P_i(x,y) = \cases{
1, &if $y=s_ix$, \cr
0, &otherwise, \cr}$$ 
with the distribution $\pi$ as given in (1.1).
Recall that the chain $K_i$ can be interpreted as follows:

{\narrower{\noindent
From $w$, try to multiply by $s_i$.  If this increases the number of 
inversions of $w$, carry out the multiplication.  If it decreases the 
the number of inversion flip a $\theta$-coin
and carry out the multiplication if the coin comes up heads. Otherwise
stay at $w$. }

}

\smallskip\noindent
The long systematic scan Metropolis chain is the chain given by
$$
K=(K_1 K_2\cdots  K_n K_n\cdots  K_2 K_1) 
\cdots
(K_1 K_2 K_2 K_1)( K_1 K_1).
$$
The following theorem bounds the rate of convergence of this Markov chain.
It shows that a single scan suffices to be close to stationarity.

\prop  Let $K$ be the long systematic scan Metropolis walk on the 
symmetric group $S_n$ defined by (4.2).  
Let $d_\lambda$, $t_\lambda$ and $c_\lambda$ be the constants given in (7.1) 
Then, for $0<\theta\le 1$,
\medskip
\item{(a)} $\displaystyle{ 
\Vert K_1^\ell/\pi-1\Vert_{{}_2}^2
= \sum_{\lambda\ne (n)}
t_\lambda d_\lambda\theta^{2\ell\left({n\choose 2}-c_\lambda\right)},
}$
\smallskip
\item{} and, with $\ell=1$,
$$\Vert K_1^1-\pi\Vert^2_{{}_{TV}}\le 
\left(e^{n^2\theta^{n/2}}-1\right)
+n!\theta^{n^2/8+5n/4},$$
which, when $\theta<1$, approaches $0$ as $n\to \infty$.
\medskip
\item{(b)} $\displaystyle{
\sum_{x\in W} \pi(x) \Vert K_x^\ell/\pi-1\Vert_{{}_2}^2 
= \sum_{\lambda\ne (n)} d_\lambda^2
\theta^{2\ell\left({n\choose 2}-c_\lambda\right)}.
}$
\smallskip
\item{} and, with $\ell=1$,
$$\sum_{x\in W} \pi(x) \Vert K_x^\ell/\pi-1\Vert_{{}_{TV}}^2 
\le 
\left(e^{n^2\theta^{n}}-1\right)
+n!\theta^{n^2/2+n},$$
which, when $\theta<1$, approaches $0$ as $n\to \infty$.
\pf
By Theorem 4.3 this walk is equivalent 
to the walk on the Iwahori-Hecke algebra $H$
defined by multiplication by $\tilde T_{w_0}^2$ with respect to the basis
$\{ \tilde T_w\ |\ w\in S_n\}$.
Thus the equalities in (a) and (b) are consequences of Theorem 4.10.

Fix $\ell=1$.
If $\lambda_1=n-j$ and $j\le n/2$ then the bound on $c_\lambda$ from Lemma 7.2(c)
gives
$$
\theta^{{\lambda_1\choose2}-{n\choose2}+2{n\choose2}-2c_\lambda} 
\le 
\theta^{j(n-j/2-1/2)-j(j-3)} 
=\theta^{j(n-3j/2+5/2)} 
\le \theta^{j(n/4+5/2)},
$$
and, by using the bounds in Lemma 7.2(a-b), it follows that 
$$
\sum_{j=1}^{n/2}\sum_{\lambda\ne (n)\atop \lambda_1 = n-j}
t_\lambda d_\lambda\theta^{2\ell\left({n\choose 2}-c_\lambda\right)}
\le \sum_{j=1}^{n/2} {n^{2j}\over j!}
\theta^{j(n/4+5/2)} 
\le e^{n^2\theta^{n/4}}-1.
$$
When $\lambda_1\le n/2$, the bound in Lemma 7.2(c) gives
$$\sum_{\lambda\atop \lambda_1\le n/2} d_\lambda^2 
\theta^{{\lambda_1\choose 2}-{n\choose2}+2{n\choose2}-2c_\lambda}
\le n!\theta^{n^2/8+5n/4}.$$
The upper bound on $\Vert K_1^1-\pi\Vert_{{}_{TV}}$ follows by combining these
expressions.  The upper bound in (b) is proved similarly.
\endpf

\subsection Short systematic scan

We now analyze the convergence of the short systematic
scan and prove Theorem 1.4 of the introduction.  
The short systematic scan Metropolis chain on the symmetric group is given by
$$K=K_1 K_2\cdots K_n K_n\cdots  K_2 K_1,$$
where $K_i$ is as in (7.3).
The Theorem shows that order $n$ short systematic scans are necessary and suffice 
to reach stationarity when starting from the identity.  In part (b${}'$) it is 
shown that for typical starting values this chain converges in order 
$\log n$ scans.

\thm  Let $K$ be the short systematic scan Metropolis algorithm
on the symmetric group defined by (4.2).  
Let $d_\lambda$, $t_\lambda$ and $c_\lambda$ be the constants given in (7.1).
Then
\medskip
\item{(a)} $\displaystyle{
\Vert K_1^\ell/\pi -1\Vert_{{}_2}^2
= \sum_{\lambda\ne (n)}
t_\lambda \sum_S \theta^{2\ell(n-1-c(S(n)))},
}$ where
\smallskip
\item{}
the sum is over standard tableaux of shape $\lambda$
and $S(n)$ denotes the box of $S$ containing $n$.
\smallskip
\item{(a${}'$)}
For $\ell = n/2-(\log n/\log \theta) + c$ with $c>0$,
$$\Vert K_1^\ell-\pi\Vert_{{}_{TV}}^2 \le \left( e^{\theta^{2c+1}}-1\right)
+n!\theta^{n^2/8 - n(\log n)/(\log \theta) +n(c+1/4)}.$$
Conversely, if $\ell< n/4$  then, for fixed $0<\theta<1$,
$\Vert K_1^\ell-\pi\Vert_{{}_{TV}}$ tends to $1$
as $n\to \infty$.
\medskip
\item{(b)} $\displaystyle{
\sum_{x\in W} \pi(x) \Vert K_1^\ell/\pi-1\Vert_{{}_2}^2
= \sum_{\lambda\ne (n)} d_\lambda \sum_S \theta^{2\ell(n-1-c(S(n)))},
}$ where
\smallskip
\item{}
the sum is over standard tableaux of shape $\lambda$
and $S(n)$ denotes the box of $S$ containing $n$.
\smallskip
\item{(b${}'$)}
For $0<\theta<1$ and $\ell = -(\log n)/(\log \theta)+c$ with $c>0$,
$$\sum_{x\in S_n} \pi(x) \Vert K_x^\ell/\pi-1\Vert_{{}_2}^2
\le \left( e^{\theta^{2c}}-1\right)+
\left({\theta^c\over e}\right)^{n}e^{1/12}\sqrt{2\pi n}.$$
\pf
By Theorem 4.3 the Markov chain $K$ is the random walk on $\{\tilde T_w\ |\ w\in S_n\}$
defined by multiplication by 
$\tilde T_{n-1}\cdots \tilde T_2\tilde T_1^2\tilde T_2\cdots \tilde T_{n-1}$
in the Iwahori-Hecke algebra $H$ corresponding to the symmetric group $S_n$.
Let $H'$ be the Iwahori-Hecke algebra corresponding to the symmetric group
$S_{n-1}$ and let $H$ be the Iwahori-Hecke algebra corresponding to 
$S_n$.  Let $w_0'$ be the longest element of $S_{n-1}$ and let $w_0$ be
the longest element of $S_n$.  The inclusion $S_{n-1}\subseteq S_n$
induces an inclusion $H'\subseteq H$ of the corresponding Iwahori Hecke
algebras.
In $H$ the generators $T_i$ are invertible with $T_i^{-1} = q^{-1}T_i+(1-q^{-1})$.
Then
$$\tilde T_{w_0}^2\tilde T_{w_0'}^{-2} 
= \tilde T_{n-1}\cdots \tilde T_2\tilde T_1^2\tilde T_2\cdots \tilde T_{n-1}$$
and so it follows from Proposition 4.9 that, in the representation
$\rho^\lambda$ of $H$ indexed by the partition $\lambda$,  the element
$$K = \tilde T_{n-1}\cdots \tilde T_2\tilde T_1^2\tilde T_2\cdots \tilde T_{n-1}
\qquad\hbox{has eigenvalues}\qquad
\theta^{n-1-c(S(n))},$$
where $S$ runs over standard tableaux of shape $\lambda$
and $S(n)$ denotes the box of $S$ containing $n$.
This determines the eigenvalues of $K^{2\ell}$ in 
the representation $\rho^\lambda$ and so 
$$\chi^\lambda_{{}_H}(K^{2\ell}) = \sum_S \theta^{2\ell(n-1-c(S(n)))},$$
where $\chi^\lambda_{{}_H}$ is the irreducible character of the Iwahori-Hecke 
algebra corresponding to the partition $\lambda$.
Parts (a) and (b) now follow from Proposition 4.8.

Using $\theta^{2\ell(n-1-c(S(n)))} \le \theta^{2\ell(n-\lambda_1)}$,
and the bound for $t_\lambda$ in Lemma 7.2 gives
$$
\Vert K_1^\ell/\pi-1\Vert_{{}_2}^2
\le 
\sum_{\lambda\ne (n)} 
\theta^{{\lambda_1\choose 2}-{n\choose2}}d_\lambda
\sum_{S} \theta^{2\ell(n-\lambda_1)} 
\le  
\sum_{\lambda\ne (n)} 
\theta^{{1\over2}(n-\lambda_1)(n-\lambda_1+4\ell+1-2n)}
d_\lambda^2, 
$$
since $d_\lambda$ is the number of standard tableaux of shape $\lambda$.
Fix $\ell = n/2-(\log n)/(\log \theta)+c$.
Then, using the bound on the sum of $d_\lambda^2$
from Lemma 7.2,
$$
\sum_{j=1}^{n/2}\sum_{\lambda\ne (n)\atop \lambda_1=n-j} 
\theta^{{1\over 2}j(j+4\ell+1-2n)} d_\lambda^2 
\le \sum_{j=1}^{n/2} 
{\theta^{{1\over 2}j(1+4\ell+1-2n+4(\log n/\log \theta))} \over j!}
=\sum_{j=1}^{\infty} {\theta^{(2c+1)j}\over j!} 
=e^{\theta^{2c+1}}-1.$$
The function
$(n-\lambda_1) (n-\lambda_1+4\ell+1-2n)$
has a minimum at $n-\lambda_1=(-1/2)(4\ell+1-2n)$.
At this minimum
$\lambda_1 = n-2(\log n)/(\log \theta)+2c+1/2\ge n/2$, and so
$$
\sum_{\lambda\atop \lambda_1\le n/2} 
\theta^{{1\over 2}(n-\lambda_1)(n-\lambda_1+4\ell+1-2n)}
d_\lambda^2 
\le \theta^{{1\over 2}(n/2)(n/2+4\ell+1-2n)} n! 
= n!\theta^{n^2/8-(n \log n)/(\log \theta)+n(c+1/4)}. 
$$
Combining these sums establishes the upper bound in (a${}'$).

To prove the lower bound in (a${}'$) let  
$$A = \{ w\in S_n \ |\ \ell(w) > 2\ell(n-1) \}
=\left\{ w\in S_n \ \Big|\ \Big|\ell(w) - {n\choose 2}\Big| < 2(n-1)(n/4-\ell) 
\right\}.$$
Since each pass of the systematic scan
can change the length of a permutation by at most $2(n-1)$,
$$K_1^\ell(A) = 0.\formula$$
From equation (2.13),
$$
E_\pi(\ell(w)) = -{(n-1)\over 1-\theta} + \sum_{j=2}^n {j\over 1-\theta^j}
= \sum_{j=2}^n j 
-{(n-1)\over 1-\theta} 
+\sum_{j=2}^n {j\theta^j\over 1-\theta^j}
= {n\choose 2}+O(n), $$
and
$$\Var_{\pi}(\ell(w)) = {(n-1)\theta\over (1-\theta)^2}
- \sum_{j=2}^n {j^2\theta^j\over (1-\theta^j)^2}
={(n-1)\theta\over (1-\theta)^2} + O(1).$$
Thus Chebychev's inequality implies that, when $n$ is large,
$$\pi(A) \sim  \left(1-{\Var_{\pi}(\ell(w))\over 
(2(n-1)(n/4-\ell))^2 }\right) \sim \left(1
-{\theta\over 4(n-1)(1-\theta)^2(n/4-\ell)^2}\right).\formula$$
If $\ell < n/4$ then the right hand side approaches
$1$ as $n\to \infty$.
Thus (7.5) and (7.6) imply that 
$\Vert K_1^\ell-\pi\Vert_{{}_{TV}} \to 1$ as $n\to \infty$ and this proves the 
second statement of (a${}'$).

\smallskip
(b${}'$)
Using the bound
$\theta^{2\ell(n-1-c(S(n)))} \le \theta^{2\ell(n-\lambda_1)},$
gives
$$
\sum_{x\in S_n} \pi(x)\Vert K_1^\ell/\pi-1\Vert_{{}_2}^2
\le
\sum_{\lambda\ne (n)} \theta^{2\ell(n-\lambda_1)} d_\lambda^2.
$$
Fix $\ell = -(\log n)/(\log \theta)+c$.
Using the bound in Lemma 7.2(a) gives
$$
\sum_{j=1}^{n/2}\sum_{\lambda\ne (n)\atop \lambda_1=n-j} 
\theta^{2\ell j} d_\lambda^2
\le 
\sum_{j=1}^{n/2} 
\theta^{2\ell j} {n^{2j}\over j!} 
=
\sum_{j=1}^{n/2} {\theta^{2j(\ell+(\log n/\log \theta))}\over j!} 
=e^{\theta^{2c}}-1, $$
and, by using the bound in Lemma 7.2(b) and the bound on $n!$ given in 
[Fe, (9.15)],
$$
\sum_{\lambda\atop \lambda_1\le n/2} 
\theta^{2\ell(n-\lambda_1)}
d_\lambda^2 
\le
\theta^{\ell n}
n! 
=
n!n^{-n}\theta^{cn}\le \left({\theta^c\over e}\right)^{n}e^{1/12}\sqrt{2\pi n}.
$$
The result follows by combining the bounds for these two sums.
\endpf

\bigskip\bigskip

\centerline{\smallcaps References}

\bigskip

\medskip
\item{[Am1]} {\smallcaps Y.\ Amit}, 
{\it  Convergence properties of the Gibbs sampler for perturbations of
Gaussians},  Ann.\ Statistics\ {\bf 24} (1995), 122--140.

\medskip
\item{[Am2]} {\smallcaps Y.\ Amit}, 
{\it  On rates of convergence of stochastic relaxation for Gaussian and
non-Gaussian distributions},  J.\ Multiv.\ Anal.\ {\bf 38} (1991), 82--99.

\medskip
\item{[AG]} {\smallcaps Y.\ Amit and U.\ Grenander}, 
{\it  Comparing sweep strategies for stochastic relaxation},  J.\ Multiv.\ Anal.\
{\bf 37} (1991), 197--222.

\medskip
\item{[AK]} {\smallcaps S.\ Ariki and K.\ Koike}, 
{\it A Hecke algebra of $(\ZZ/r\ZZ)\wr S_n$ and construction of its irreducible
representations}, Adv.\ Math.\ {\bf 106} (1994), 216--243.

\medskip
\item{[Be]} {\smallcaps E.\ Belsley}, 
{\it Rates of convergence of random walk on distance regular graphs},
Probab.\ Theory Related Fields {\bf 112} (1998), 493--533. 

\medskip
\item{[Be2]} {\smallcaps E.\ Belsley},
{\sl Rates of convergence of Markov chains related to association schemes},
Ph.D.\ thesis, Harvard, 1993.

\medskip
\item{[BF]} {\smallcaps P.\ Barone and A.\ Frigessi}, {\it Improving stochastic
relaxation for Gaussian random fields}, Probab. Eng. Inform. Sci. {\bf 4}
(1990), 369--389.

\medskip
\item{[BHV]} {\smallcaps L.\ Billera, S.\ Holmes, K. Vogtman},
{\it Geometry of the space of phylogenetic trees}, Technical Report,
Department of Statistics, Stanford University, 1999.

\medskip
\item{[BS]} {\smallcaps E.\ Brieskorn and K.\ Saito}, {\it Artin-Gruppen und Coxeter-Gruppen},
Invent.\ Math.\ {\bf 17} (1972), 245--271. 

\medskip
\item{[Bou]} {\smallcaps N.\ Bourbaki}, {\sl Groupes et alg\`ebres de Lie, Ch. 4,5 et 6},
Masson, Paris, 1981.

\medskip
\item{[Bw]} {\smallcaps K.\ Brown}, {\sl Buildings}, Springer-Verlag, New York-Berlin,
1989.

\medskip
\item{[Ca]} {\smallcaps R.\ Carter}, {\sl Finite groups of Lie type -- 
Conjugacy classes and complex characters},
Wiley-Interscience, John Wiley \& Sons, Inc., New York, 1985.

\medskip
\item{[Cr]} {\smallcaps D.\ Critchlow}, {\sl Metric methods for analyzing
partially ranked data}, Lecture Notes in Statistics {\bf 34}, Springer-Verlag,
Berlin, 1985.  

\medskip
\item{[CR]} {\smallcaps C.\ Curtis and I.\ Reiner}, {\sl Methods of 
representation theory -- with applications to finite groups and orders,
Volumes I and II}, Wiley-Interscience, John Wiley \& Sons, Inc., 
New York, 1981 and 1987.

\medskip
\item{[De]} {\smallcaps P.\ Deligne}, {\it Les immeubles des groupes de 
tresses g\'en\'eralis\'es}, Invent.\ Math.\ {\bf 17} (1972), 273--302. 

\medskip
\item{[D]} {\smallcaps P.\ Diaconis}, {\sl Group representations in probability and
statistics}, Institute of Mathematical Statistics, Hayward, CA, 1988.

\medskip
\item{[DH]} {\smallcaps P.\ Diaconis and P.\ Hanlon},
{\it Eigen-analysis for some examples of the Metropolis algorithm}, in 
{\sl Hypergeometric functions on domains of positivity, Jack polynomials, 
and applications} (Tampa, FL, 1991), 99--117, Contemp.\ Math., {\bf 138}, 
Amer.\ Math.\ Soc., Providence, RI, 1992. 

\medskip
\item{[Det]} {\smallcaps P.\ Diaconis, S.\ Holmes, S.\ Janson, S.\ Lalley,
R.\ Pemantle}, {\it Metrics on compositions and coincidences among renewal sequences},
in {\sl Random discrete structures}, D.\ Aldous and R. Pemantle eds., Springer,
New York, 1996.

\medskip
\item{[DS]} {\smallcaps P.\ Diaconis and L.\ Saloff-Coste},
{\it What do we know about the Metropolis algorithm?},
J.\ Computer and System Sciences {\bf 57} (1998), 20-36.

\medskip
\item{[Fe]} {\smallcaps W.\ Feller},
{\sl An introduction to probability theory and its applications},
J.\ Wiley and Sons, New York, 1968.

\medskip
\item{[Fi]} {\smallcaps G.\ Fishman}, 
{\it  Coordinate selection rules for Gibbs sampling},  Ann.\ Appl.\ Probab.\ 
{\bf 6} (1996) 444--465.

\medskip
\item{[FV]} {\smallcaps M.\ Figner and J. Verducci},
{\it Probability models and statistical analyses for ranking data},
Springer Lecture Notes in Statistics {\bf 80}, Springer, New York, 1993.

\medskip
\item{[Fu]} {\smallcaps W.\ Fulton},
{\sl Young tableaux, With applications to representation theory 
and geometry}, London Math.\ Soc.\ Student Texts {\bf 35}, 
Cambridge Univ.\ Press, Cambridge, 1997.

\medskip
\item{[GS]} {\smallcaps J.\ Goodman and A.\ Sokal}, 
{\it Multigrid Monte-Carlo method.  Conceptual foundations},
Phys.\ Rev.\ D {\bf 40} (1989), 2035--2071.

\medskip
\item{[HH]} {\smallcaps J.\ Hammersley and D.\ Handscomb}, {\sl Monte Carlo Methods},
Chapman and Hall, London, 1964.

\medskip
\item{[Hf]}  {\smallcaps P.N.\ Hoefsmit}, {\sl Representations of Hecke algebras
of finite groups with BN-pairs of classical type}, Ph.D. Thesis, University 
of British Columbia, 1974.

\medskip
\item{[Hu]} {\smallcaps J.\ Humphreys},
{\sl Reflection groups and Coxeter groups}, Cambridge University Press,
1990.

\medskip
\item{[KS]} {\smallcaps J.\ Kemeny and L.\ Snell},
{\sl Finite Markov chains}, The University Series in Undergraduate Mathematics,
Van Nostrand, Princeton, 1960.

\medskip
\item{[KSo]} {\smallcaps R.\ Kilmoyer and L.\ Solomon},
{\it On the theorem of Feit-Higman}, 
J.\ Combinatorial Theory A {\bf 15} (1973), 310-322.

\medskip
\item{[Le1]} {\smallcaps G.\ Letac}, {\it Probl\`emes classiques de probabilit\'e
sur un couple de Gelfand}, in Lect. Notes in Math. {\bf 861}, Springer-Verlag,
New York, 1981.

\medskip
\item{[Le2]} {\smallcaps G.\ Letac}, {\it Les fonctions sph\'eriques d'un couple de 
Gelfand sym\'etrique et les cha\^ines de Markov}, Advances in Appl.\ Prob.\ {\bf 14}
(1982), 272-294.

\medskip
\item{[Ma]} {\smallcaps C.L.\ Mallows}, {\it Non-null ranking models I},
Biometrika {\bf 44} (1957), 114-130.

\medskip
\item{[Mac]}  {\smallcaps I.G.\ Macdonald}, {\sl Symmetric functions and Hall polynomials},
Second edition, Oxford Mathematical Monographs, Oxford University Press, New York, 1995.

\medskip
\item{[Mar]} {\smallcaps J. Marden}, {\sl Analyzing and modeling rank data}, 
Monographs on Statistics and Applied Probability {\bf 64}, Chapman \& Hall,
London, 1995.

\medskip
\item{[MR]} {\smallcaps N.\ Metropolis, A.\ Rosenbluth, M.\ Rosenbluth, A.\ Teller
and E.\ Teller}, {\it Equations of state calculations by fast computing machines},
J.\ Chem.\ Phys.\ {\bf 21} (1953), 1087-1092.

\medskip
\item{[Pa]} {\smallcaps I. Pak}, {\it When and how $n$ choose $k$},
in  {\sl Randomization methods in algorithm design (Princeton, NJ, 1997)}, 191--238,
DIMACS Ser. Discrete Math. Theoret. Comput. Sci. {\bf 43}, Amer. Math. Soc., 
Providence, RI, 1999.

\medskip
\item{[Ra]} {\smallcaps A.\ Ram}, {\it Seminormal representations of Weyl groups 
and Iwahori-Hecke algebras}, Proc.\ London Math.\ Soc.\ {\it (3)} {\bf 75} (1997),
99--133.   

\medskip
\item{[RS]} {\smallcaps G.\ Roberts and S.\ Sahk}, 
{\it  Updating schemes, correlation schemes, correlation structure, blocking
and parameterization for the Gibb sampler},  J.\ Royal Statist.\ Soc.\ {\bf 59}
(1997), 291--318.

\medskip
\item{[RX]} {\smallcaps K.\ Ross and D.\ Xu},
{\it Hypergroup deformations and Markov chains}, J. Theor.\ Probability {\bf 7}
(1994), 813--830.

\medskip
\item{[SB]} {\smallcaps W.\ Shannon and D.\ Banks}, 
{\it Combining classification trees using MLE}, Statistics in Medicine {\bf 18} (1999),
727-740.

\medskip
\item{[SC]} {\smallcaps L.\ Saloff-Coste},
{\sl Lectures on probability theory and statistics} (Saint-Flour, 1996), 301--413, 
Lecture Notes in Math. {\bf 1665}, Springer, Berlin, 1997.

\medskip
\item{[Sc]} {\smallcaps C.\ Schoolfield},
{\sl Random walks on wreath products of groups and Markov chains on related
homogeneous spaces}, Ph.D.\ Thesis, 
Dept.\ of Appl.\ Math., Johns Hopkins Univ., 1998.

\medskip
\item{[Si]} {\smallcaps J.\ Silver}, 
{\sl Weighted Poincar\'e and exhaustive approximation techniques for scaled
Metropolis-Hastings and spectral total variation convergence bounds in infinite
commutable Markov chain theory},
Ph.D.\ Thesis, Dept.\ of Math., Harvard University, 1996.

\vfill\eject
\end

The formula for $t_\lambda$ for the hypercube can be checked as follows.
$$\eqalign{
\sum_{x\in (\ZZ/2\ZZ)^n} q^{-\ell(x)}\chi_{{}_H}^\lambda(T_{x^{-1}})\chi_{{}_H}^\lambda(T_x)
&= \sum_{x\in (\ZZ/2\ZZ)^n} q^{-|x|}\chi_{{}_H}^\lambda(T_x)^2 \cr
&= \sum_{x\in (\ZZ/2\ZZ)^n} q^{-|x|}q^{2(|x|-\lambda\cdot x)}(-1)^{2\lambda\cdot x} \cr
&= \sum_{x\in (\ZZ/2\ZZ)^n} q^{|x|-2\lambda\cdot x} \cr
&=\prod_{i=1}^n (q^{1-2\lambda_i}-1) \cr
&=\prod_{\lambda_i\ne 0} (q^{-1}+1)\prod_{\lambda_i=0} (q+1)
=(1+q)^n q^{-|\lambda|}. \cr
}$$
So
$$t_\lambda = {(1+q)^n\over (1+q)^nq^{-|\lambda|}} = q^{|\lambda|}.$$

\section 8.  The hyperoctahedral group

While we will not carry out the details
we here present the needed tools for proving theorems as above for the hyperoctahedral
group.

The hyperoctahedral group $WB_n$ can also be studied using
present techniques.  Ordinary (undeformed) walks on $WB_n$ arise in
the analysis of card shuffling schemes where the order and the
orientation of the back designs is of interest.  They also arise in 
DNA applications.  C. Schofield [Sc??] has carried out a careful analysis
which should also be applicable to deformed walks.  In the present
brief section we give the ingredients for such an analyasis.

The {\it hyperoctahedral group}, or {\it Weyl group of type $B_n$}, is the 
group $WB_n$ of signed permutations, i.e. permutations of 
$\{-n,-(n-1),\ldots, -2,-1,1,2,\ldots,n-1,n\}$ such that 
$w(-i)=-w(i)$.  The group $WB_n$ is generated by
$$s_1=(-1,1)
\quad\hbox{and}\quad s_i=(i,i+1)(-i,-(i+1)),
\quad 2\le i\le n,$$
and
$w_0=\pmatrix{1 &2 &\cdots &n-1 &n\cr -1 &-2 &\cdots &-(n-1) &-n \cr}$
is the longest element of $S_n$ with length $\ell(w_0)=n^2$.

The irreducible representations of the Iwahori-Hecke algebra are indexed
by pairs of partitions $\lambda=(\alpha,\beta)$ with $n$ boxes total.
The {\it sign}, {\it content}, {\it hook length} and {\it split hook length}
of a box $b$ in $\lambda$ are given, respectively, by
$$\eqalign{
{\rm sgn}(b) &= 
\cases{ 1, &if $b$ is in $\alpha$, \cr
       -1, &if $b$ is in $\beta$, \cr}
\cr
\cr
c(b) &= j-i, \qquad \hbox{if $b$ is in position $(i,j)$ (in either $\alpha$ or $\beta$)}, \cr
\cr
h_b &= \cases{
\alpha_i-j+\alpha_j'-i+1, &if $b$ is in position $(i,j)$ of $\alpha$, \cr
\cr
\beta_i-j+\beta_j'-i+1, &if $b$ is in position $(i,j)$ of $\beta$, \cr
} \cr
\cr
\tilde h_b &= \cases{
\alpha_i-j+\beta_j'-i+1, &if $b$ is in position $(i,j)$ of $\alpha$, \cr
\cr
\beta_i-j+\alpha_j'-i+1, &if $b$ is in position $(i,j)$ of $\beta$. \cr
} \cr
}$$
For each positive integer
$k$ let $[k] = (q^k-1)/(q-1)$ and $[2k]!!=[2k][2k-2]\cdots [4][2]$.
With these notations
$$c_\lambda = \sum_{b\in \lambda} {\rm sgn}(b) +2c(b)
\qquad\hbox{and}\qquad
t_\lambda={[2n-1]!!\over \displaystyle{\prod_{b\in\lambda} [h_b][\tilde h_b]} }\ ,$$
see [Ra, Prop. 4.8] and [Hf, Cor. 3.4.11].